\documentclass[12pt,a4paper,reqno]{amsart}

\usepackage[T1]{fontenc}
\usepackage[utf8]{inputenc}
\usepackage{lmodern}

\usepackage{makecell}
\usepackage[12pt]{moresize}

\usepackage[
    margin=1.2in
]{geometry}

\usepackage{amsmath,amsthm}
\usepackage{amssymb}
\usepackage{mathptmx}
\IfFileExists{omlntxmi.fd}{%
  \DeclareSymbolFont{newtxletters}{OML}{ntxmi}{m}{it}%
}{%
  \DeclareSymbolFont{newtxletters}{OML}{cmm}{m}{it}%
}
\DeclareMathSymbol{\partial}{\mathord}{newtxletters}{"40}
\DeclareMathAlphabet{\mathcal}{OMS}{cmsy}{m}{n}
\DeclareSymbolFont{cmletters}{OML}{cmm}{m}{it}
\DeclareMathSymbol{v}{\mathalpha}{cmletters}{`v}
\DeclareSymbolFont{cmoperators}{OT1}{cmr}{m}{n}
\DeclareMathSymbol{\Gamma}{\mathalpha}{cmoperators}{"00}
\DeclareMathSymbol{\Delta}{\mathalpha}{cmoperators}{"01}
\DeclareMathSymbol{\Theta}{\mathalpha}{cmoperators}{"02}
\DeclareMathSymbol{\Lambda}{\mathalpha}{cmoperators}{"03}
\DeclareMathSymbol{\Xi}{\mathalpha}{cmoperators}{"04}
\DeclareMathSymbol{\Pi}{\mathalpha}{cmoperators}{"05}
\DeclareMathSymbol{\Sigma}{\mathalpha}{cmoperators}{"06}
\DeclareMathSymbol{\Upsilon}{\mathalpha}{cmoperators}{"07}
\DeclareMathSymbol{\Phi}{\mathalpha}{cmoperators}{"08}
\DeclareMathSymbol{\Psi}{\mathalpha}{cmoperators}{"09}
\DeclareMathSymbol{\Omega}{\mathalpha}{cmoperators}{"0A}
\usepackage{graphicx}

\usepackage{
    color,
    tikz,
    xfrac,
    etoolbox,
    xparse,
    pgffor,
    longtable,
    booktabs,
    makecell,
    array,
    placeins,
    xcolor,
    microtype,
    multirow,
    tabularx
}

\IfFileExists{siunitx.sty}{\usepackage{siunitx}}{}
\IfFileExists{dsfont.sty}{\usepackage{dsfont}}{}
\IfFileExists{bropd.sty}{\usepackage{bropd}}{}
\usepackage{needspace}

\usepackage[shortlabels]{enumitem}

\setlist{
    font=\upshape,
    itemsep=0.25\baselineskip,
    leftmargin=\parindent,
    parsep=0pt,
    topsep=0.25\baselineskip
}

\usetikzlibrary{
    arrows.meta,
    calc,
    decorations.pathreplacing,
    positioning
}

\definecolor{freeviolet}{RGB}{0,0,255}
\definecolor{contactred}{RGB}{0,0,0}
\definecolor{phasegreen}{RGB}{255,255,255}

\newtheoremstyle{plain}
    {6.5pt}
    {6.5pt}
    {\itshape}
    {}
    {\bfseries}
    {.}
    {.5em}
    {}

\newtheoremstyle{definition}
    {6.5pt}
    {6.5pt}
    {\normalfont}
    {}
    {\bfseries}
    {.}
    {.5em}
    {}

\theoremstyle{plain}

\newtheorem{theorem}{Theorem}[section]
\newtheorem{lemma}[theorem]{Lemma}

\newtheorem{proposition}[theorem]{Proposition}

\theoremstyle{definition}

\newtheorem{definition}[theorem]{Definition}
\newtheorem{example}[theorem]{Example}
\newtheorem{remark}[theorem]{Remark}

\theoremstyle{plain}

\newtheorem{innercustomgeneric}{\customgenericname}
\providecommand{\customgenericname}{}

\newcommand{\newcustomtheorem}[2]{%
    \newenvironment{#1}[1]
    {%
        \renewcommand\customgenericname{#2}%
        \renewcommand\theinnercustomgeneric{##1}%
        \innercustomgeneric
    }
    {\endinnercustomgeneric}
}

\newcustomtheorem{maintheorem}{Theorem}

\numberwithin{equation}{section}

\usepackage[
    pdfa,
    colorlinks,
    allcolors=blue
]{hyperref}

\makeatletter

\let\keywords\relax
\let\msc\relax
\let\correspondingauthor\relax
\let\abstract\relax

\gdef\@keywords{}
\gdef\@msc{}
\gdef\@correspondingauthor{}

\newcommand{\keywords}[1]{%
    \gdef\@keywords{%
        \textit{Keywords:} \ignorespaces#1\unskip.%
    }%
}

\newcommand{\msc}[2][]{%
    \gdef\@msc{%
        \textit{Mathematics Subject Classification 2020:} \ignorespaces#2\unskip%
        \ifx&#1&\else\ (#1)\fi.%
    }%
}

\newcommand{\correspondingauthor}[2]{%
    \gdef\@correspondingauthor{%
        \textit{Corresponding author:} #1,
        \href{mailto:#2}{#2}.%
    }%
}

\newcommand{\ems@titlefootnote}[1]{%
    \insert\footins{%
        \reset@font
        \footnotesize
        \interlinepenalty\interfootnotelinepenalty
        \splittopskip\footnotesep
        \splitmaxdepth\dp\strutbox
        \floatingpenalty\@MM
        \hsize\columnwidth
        \@parboxrestore
        \color@begingroup
            \noindent
            \strut
            #1\par
        \color@endgroup
    }%
}

\newbox\emsabstractbox

\newenvironment{abstract}{%
    \global\setbox\emsabstractbox\vbox\bgroup
        \normalcolor
        \small
        \noindent
        \strut
        \textbf{Abstract.}\enskip
        \ignorespaces
}{%
        \unskip
        \strut
        \par
    \egroup
}

\renewcommand{\maketitle}{%

    \thispagestyle{empty}

    \begingroup

        \vspace*{26pt}

        \centering

        {\Large\bfseries
            \@title
            \par
        }

        \vspace{8mm}

        {\large
            \@author
            \par
        }

        \vspace{5mm}

    \endgroup

    \ifx\@msc\@empty
    \else
        \ems@titlefootnote{\@msc}%
    \fi

    \ifx\@keywords\@empty
    \else
        \ems@titlefootnote{\@keywords}%
    \fi

    \ifx\@correspondingauthor\@empty
    \else
        \ems@titlefootnote{\@correspondingauthor}%
    \fi

    \ifvoid\emsabstractbox
    \else
        \noindent
        \unvbox\emsabstractbox
        \par
    \fi

    \vspace{26pt}
}

\newcounter{ems@authorcount}
\gdef\@author{}
\gdef\ems@authorids{}

\newcommand{\Author}[3]{%
    \stepcounter{ems@authorcount}%
    \listgadd{\ems@authorids}{#1}%
    \begingroup
        \def\givenname##1{\csgdef{auth@#1@first}{##1}}%
        \def\surname##1{\csgdef{auth@#1@last}{##1}}%
        \def\mrid##1{%
            \ifblank{##1}{}{\csgdef{auth@#1@mrid}{##1}}%
        }%
        \def\zblid##1{%
            \ifblank{##1}{}{\csgdef{auth@#1@zblid}{##1}}%
        }%
        \def\orcid##1{%
            \ifblank{##1}{}{\csgdef{auth@#1@orcid}{##1}}%
        }%
        #2%
    \endgroup
    \ifnum\value{ems@authorcount}=1
        \gdef\@author{%
            \csuse{auth@#1@first}\space\csuse{auth@#1@last}%
        }%
    \else
        \g@addto@macro\@author{%
            ,\space\csuse{auth@#1@first}\space\csuse{auth@#1@last}%
        }%
    \fi
}

\def\ems@setmeta#1#2#3{%
    \ifblank{#3}{}{%
        \csgdef{affil@#1@#2}{#3}%
        \csgdef{affil@#1@present}{1}%
    }%
}

\def\ems@makeaffilmacro#1{%
    \expandafter\def\csname #1\endcsname##1{%
        \@ifnextchar\bgroup
            {\csuse{ems@#1@two}{##1}}%
            {\csuse{ems@#1@one}{##1}}%
    }%
    \expandafter\def\csname ems@#1@two\endcsname##1##2{%
        \ems@setmeta{\ems@currentaffil@id @##1}{#1}{##2}%
    }%
    \expandafter\def\csname ems@#1@one\endcsname##1{%
        \ems@setmeta{\ems@currentaffil@id @1}{#1}{##1}%
    }%
}

\newcommand{\Affil}[2]{%
    \begingroup
        \def\ems@currentaffil@id{#1}%
        \ems@makeaffilmacro{department}%
        \ems@makeaffilmacro{organisation}%
        \ems@makeaffilmacro{address}%
        \ems@makeaffilmacro{post}%
        \ems@makeaffilmacro{city}%
        \ems@makeaffilmacro{country}%
        \ems@makeaffilmacro{affemail}%
        #2%
    \endgroup
}

\def\ems@emitfield#1{%
    \ifcsdef{affil@\ems@affkey @#1}{%
        \ems@sep
        \csuse{affil@\ems@affkey @#1}%
        \def\ems@sep{,\space}%
    }{}%
}

\def\ems@printid#1#2#3{%
    \ifcsdef{auth@\ems@aid @#1}{%
        \edef\ems@tempid{\csuse{auth@\ems@aid @#1}}%
        \ems@idsep
        #2~\href{#3\ems@tempid}{\ems@tempid}%
        \def\ems@idsep{\quad}%
    }{}%
}
\def\ems@checkid#1{%
    \ifcsdef{auth@\ems@aid @#1}{\def\ems@hasid{1}}{}%
}

\def\ems@printauthor#1{%
    \begingroup
        \def\ems@aid{#1}%
        \Needspace{5\baselineskip}%
        \noindent
        {\bfseries\csuse{auth@#1@first}\space\csuse{auth@#1@last}}%
        \par\nobreak
        \foreach \j in {1,2,3,4}{%
            \edef\ems@affkey{\ems@aid @\j}%
            \ifcsdef{affil@\ems@affkey @present}{%
                \def\ems@sep{}%
                \ems@emitfield{department}%
                \ems@emitfield{organisation}%
                \ifx\ems@sep\@empty\else\par\nobreak\fi
                \def\ems@sep{}%
                \ems@emitfield{address}%
                \ems@emitfield{post}%
                \ems@emitfield{city}%
                \ems@emitfield{country}%
                \ifx\ems@sep\@empty\else\par\fi
                \ifcsdef{affil@\ems@affkey @affemail}{%
                    \edef\ems@tempemail{%
                        \csuse{affil@\ems@affkey @affemail}%
                    }%
                    \href{mailto:\ems@tempemail}{\ems@tempemail}\par
                }{}%
            }{}%
        }%
        \def\ems@hasid{0}%
        \forcsvlist{\ems@checkid}{zblid,mrid,orcid}%
        \ifnum\ems@hasid=1
            \def\ems@idsep{}%
            \noindent Author IDs:\space
            \ems@printid{zblid}{zbMATH}{https://zbmath.org/authors/}%
            \ems@printid{mrid}{MR}{https://mathscinet.ams.org/mathscinet/MRAuthorID/}%
            \ems@printid{orcid}{ORCID}{https://orcid.org/}%
            \par
        \fi
        \addvspace{13pt}%
    \endgroup
}

\AtEndDocument{%
    \ifnum\value{ems@authorcount}>0
        \par\addvspace{19.5pt}%
        \begingroup
            \small\raggedright
            \forlistloop{\ems@printauthor}{\ems@authorids}%
        \endgroup
    \fi
}

\makeatother

\newcommand{\dx}[1]{%
    \mathop{}\!\mathrm{d}#1
}

\newcommand{\dH}[1]{%
    \mathop{}\!\mathrm{d}\mathcal{H}^{#1}
}

\def\Xint#1{\mathchoice
 {\XXint\displaystyle\textstyle{#1}}
 {\XXint\textstyle\scriptstyle{#1}}
 {\XXint\scriptstyle\scriptscriptstyle{#1}}
 {\XXint\scriptscriptstyle\scriptscriptstyle{#1}}\!\int}
\def\XXint#1#2#3{{\setbox0=\hbox{\(#1{#2#3}{\int}\)}
 \vcenter{\hbox{\(#2#3\)}}\kern-.53\wd0}}
\def\barint{\Xint-}

\newcommand{\R}{\mathbb{R}}
\newcommand{\C}{\mathbb{C}}

\let\ge\geqslant
\let\le\leqslant

\newcommand{\mres}{%
    \mathbin{%
        \vrule
            height 1.6ex
            depth 0pt
            width 0.13ex
        \vrule
            height 0.13ex
            depth 0pt
            width 1.3ex
    }%
}

\Author{1}{%

    \givenname{Lili}
    \surname{Du}

    \mrid{}
    \zblid{}
    \orcid{}

}{L.~Du}

\Author{2}{%

    \givenname{Chunlei}
    \surname{Yang}

    \mrid{}
    \zblid{}
    \orcid{}

}{C.~Yang}

\Author{3}{%

    \givenname{Fan}
    \surname{Zhang}

    \mrid{}
    \zblid{}
    \orcid{}

}{F.~Zhang}

\Affil{1}{%

    \department{Department of Mathematics}

    \organisation{Sichuan University}

    \address{}

    \post{610065}

    \city{Chengdu}

    \country{China}

    \affemail{dulili@scu.edu.cn}
}

\Affil{2}{%

    \department{School of Mathematical Sciences}

    \organisation{Shenzhen University}

    \address{}

    \post{518060}

    \city{Shenzhen}

    \country{China}
    
    \affemail{yangchunlei@szu.edu.cn}
}

\Affil{3}{%

    \department{School of Mathematics}

    \organisation{North University of China}

    \address{}

    \post{030051}

    \city{Taiyuan}

    \country{China}
    
    \affemail{20220182@nuc.edu.cn}
}

\definecolor{freeviolet}{RGB}{0,0,0}
\DeclareMathOperator{\diver}{div}
\DeclareMathOperator{\supp}{supp}

\newcommand{\loc}{\mathrm{loc}}
\newcommand{\dz}{\partial_z}
\newcommand{\dbz}{\partial_{\bar z}}
\newcommand{\CT}{\mathcal T}
\newcommand{\BT}{\mathcal B}
\newcommand{\US}{u_{\mathrm{st}}}

\title{Optimal local gradient estimate for gravity water waves and its application to the Stokes conjecture}
\date{}
\msc{35R35, 35B65, 76B15}
\correspondingauthor{F.~Zhang}{20220182@nuc.edu.cn}

\hypersetup{
  pdftitle={Optimal local gradient estimate for gravity water waves and its application to the Stokes conjecture},
  pdfauthor={L. Du, C. Yang, F. Zhang},
  pdfsubject={Local gradient estimates from full-domain stationarity and Stokes corners},
  pdfkeywords={Stokes conjecture, irrotational gravity waves, stagnation points, domain variation identity, optimal growth}
}

\begin{document}

\begin{abstract}
\setlength{\emergencystretch}{1em}
In 1880, G. Stokes [\emph{Mathematical and Physical Papers}, Vol.~I]
studied a two-dimensional gravity water wave and conjectured that 
the wave crest forms a symmetric corner of angle \(2\pi/3\) at the stagnation point. 
The conjecture was first proved by
Amick, Fraenkel, and Toland [\emph{Acta Math.} 148 (1982)] and,
independently, by Plotnikov [\emph{Dinamika Sploshn. Sredy}
57 (1982)] via the Nekrasov integral equation. Their proofs rely
on structural assumptions, including the isolation of stagnation
points and the symmetry and monotonicity of the wave profile.
In 2011, V\u{a}rv\u{a}ruc\u{a} and Weiss
[\emph{Acta Math.} 206 (2011)] introduced a novel geometric approach
and gave a new proof of the conjecture without relying on the
structural assumptions mentioned above, under a local gradient bound for the two-dimensional stream function.

In this paper, we establish an optimal local gradient estimate for general
variational solutions of the free boundary problem.
Our approach adapts the Hopf differential method of
De Philippis, Hirsch, and Nahon [\href{https://arxiv.org/abs/2609.12180v2}{arXiv:2609.12180v2}]
to the gravitational setting, and combines it with their abstract
approximation theorem to turn the holomorphic approximation of the Hopf
differential into the local gradient bounds.
An important new observation is that the weighted Cauchy--Beurling
correction, combined with vanishing on the lower half plane, supplies the
additional height decay needed to obtain the local gradient bounds.
As a consequence, we remove the a priori growth hypothesis from the
geometric analysis in [\emph{Acta Math.} 206 (2011)]. For general variational
solutions, we establish the Stokes corner or full-density blow-up alternative.
Under the stated topological and interface hypotheses, we obtain the Stokes
corner and two free boundary branches that are \(C^1\) up to the stagnation point.
\end{abstract}

\maketitle

\enlargethispage{\baselineskip}
\tableofcontents

\section{Introduction}

The mathematical study of fluid motion dates back to Euler's work
in the 1750s. Since then, a wide range of mathematical tools and
ideas have been developed to study gravity water waves. We begin
by recalling the Stokes conjecture, one of the best-known conjectures
in the study of water waves, and reviewing the development of its
proofs and the techniques they employ.

\subsection{The water wave problem and the Stokes conjecture }

We consider an inviscid, incompressible, and irrotational two-dimensional
fluid under the influence of gravity, with a free surface separating
the fluid from the air. Let \(D(t)\subset\R^2\) be the fluid region.
The velocity field \((u,v)\) and the pressure \(P\) satisfy the
Euler equations
\begin{equation}\label{eq:euler}
	\left\{
	\begin{aligned}
		u_t+uu_x+vu_y&=-P_x,\\
		v_t+uv_x+vv_y&=-P_y-g,\\
		u_x+v_y&=0,\\
		u_y-v_x&=0,
	\end{aligned}
	\right.
	\qquad\hbox{in }D(t),
\end{equation}
where \(g>0\) is the gravitational acceleration and subscripts denote
partial derivatives. We denote by \(\partial_aD(t)\) the free part
of the boundary, which is in contact with the air. The kinematic
boundary condition states that this interface moves with the fluid,
i.e.,
\[
	V=(u,v)\cdot\nu\qquad\hbox{on }\partial_aD(t),
\]
where \(V\) is the normal velocity of the surface and \(\nu\) is its
outward unit normal. The pressure \(P\) satisfies the dynamic
boundary condition
\[
	P\quad\hbox{is locally constant on }\partial_aD(t).
\]

To study the shape of the wave near a stagnation point, we consider
a \emph{traveling wave} that moves at a constant speed \(c\) without
changing its shape. More precisely, we write
\[
	D(t)=D+ct(1,0),\qquad
	(u,v,P)(x,y,t)
	=(\widetilde u+c,\widetilde v,\widetilde P)(x-ct,y),
\]
where the domain \(D\) and the functions
\((\widetilde u,\widetilde v,\widetilde P)\) are independent of time.
In this frame, the system \eqref{eq:euler} and the boundary conditions
take the form
\[
	\left\{
	\begin{aligned}
		\widetilde u\widetilde u_x+\widetilde v\widetilde u_y
		&=-\widetilde P_x,\\
		\widetilde u\widetilde v_x+\widetilde v\widetilde v_y
		&=-\widetilde P_y-g,\\
		\widetilde u_x+\widetilde v_y&=0,\\
		\widetilde u_y-\widetilde v_x&=0,\\
		(\widetilde u,\widetilde v)\cdot\nu&=0
		\quad\hbox{on }\partial_aD,\\
		\widetilde P&\quad\hbox{is locally constant on }\partial_aD.
	\end{aligned}
	\right.
\]
Incompressibility implies the local existence of a stream function
\(\psi\) satisfying
\[
	\psi_x=-\widetilde v,\qquad
	\psi_y=\widetilde u.
\]
We work near a connected portion of the free surface, in a fluid
neighborhood where \(\psi\) is single valued and continuous up to
that portion. The kinematic condition makes \(\psi\) constant
there, and we normalize this constant to zero. We assume that
\(\psi>0\) in the fluid part of the neighborhood and extend it by
zero into the air region. Thus the fluid region is locally identified
with \(\{\psi>0\}\). In terms of \(\psi\), we have
\[
	\Delta\psi=0\qquad\hbox{in }\{\psi>0\}.
\]
Bernoulli's law yields
\[
	\widetilde P+\frac12|\nabla\psi|^2+gy
	\quad\hbox{is locally constant in }D.
\]
The dynamic boundary condition then implies the Bernoulli free boundary condition
\[
	|\nabla\psi|^2+2gy\quad\text{is locally constant on }\partial_aD.
\]
After a standard normalization of the Bernoulli constant, the height,
and the gravitational constant \(g\), we obtain the following Bernoulli-type free boundary problem:
\begin{equation}\label{eq:physical-intro}
	\left\{
	\begin{aligned}
		\Delta\psi&=0\quad &&\hbox{in }\{\psi>0\},\\
		\psi&=0 \quad &&\hbox{on }\partial\{\psi>0\},\\
		|\nabla\psi|^2&=-y\quad &&\hbox{on }\partial\{\psi>0\}.
	\end{aligned}
	\right.
\end{equation}
We are particularly interested in the behavior of solutions to
\eqref{eq:physical-intro} near stagnation points on the free boundary, namely, points
where the relative velocity \((\psi_y,-\psi_x)\) vanishes. It follows
from the Bernoulli-type free boundary condition \(|\nabla\psi|^2=-y\) that
stagnation points lie in the intersection
\(\partial\{\psi>0\}\cap\{y=0\}\). Moreover, the free surface lies
in the closed lower half plane \(\{y\le 0\}\). We assume that the
fluid also lies below this line in the neighborhood under
consideration, so the zero extension vanishes on \(\{y\ge0\}\).
The study of the
local geometry near stagnation points dates back to Stokes, who
formulated the following conjecture \cite[pp.~225--228]{St80}:
\begin{quotation}
\dots \textit{the form in which the curvature at the vertex becomes
infinite, and the vertex becomes a multiple point where the two
branches with which alone we are concerned enclose an angle of
\(120^\circ\)}\dots
\end{quotation}
In the same paper \cite{St80}, Stokes supported his conjecture with
a formal argument, approximating the stream function \(\psi\) by
a homogeneous function and thereby obtaining the \(120^\circ\)
corner at the stagnation point.

The classical approach to the conjecture is based on conformal
mappings and Nekrasov's integral equation \cite{Nek21}. The first
rigorous proofs were given independently by Amick, Fraenkel, and
Toland \cite{AFT82} and by Plotnikov \cite{P82}. Both proofs rely
on the integral equation. Moreover, some structural assumptions on
the solutions were required: the stagnation point is isolated,
and the wave profile is symmetric and monotone on either side of
the stagnation point.

\subsection{The geometric approach to the Stokes conjecture}
In their breakthrough work \cite{VW11}, V\u{a}rv\u{a}ruc\u{a} and
Weiss studied the problem \eqref{eq:physical-intro} in the original
physical variables and developed a geometric proof of the conjecture
without the structural assumptions of the earlier
proofs. We briefly recall their formulation and the main ideas of
their proof in the following. Let \(x=(x_1,x_2)\) denote the reflected position, whose components 
are related to the original physical coordinates \((x,y)\) by \(x_1=x\) and \(x_2=-y\), and define
\[
	u(x_1,x_2)=\psi(x_1,-x_2).
\]
We then obtain the following equivalent problem for \(u\):
\begin{equation}\label{eq.problem}
	\left\{
	\begin{aligned}
		\Delta u&=0\quad &&\hbox{in }\Omega\cap\{u>0\},\\
		u&=0\quad&&\text{on }\Omega\cap\partial\{u>0\},\\
		|\nabla u|^2&=x_2\quad
		&&\text{on }\Omega\cap\partial\{u>0\}.
	\end{aligned}
	\right.
\end{equation}
Here \(\Omega\subset\R^2\) is a bounded domain that intersects \(\{x_2=0\}\).
Solutions of \eqref{eq.problem} can be viewed as critical points
of the energy functional
\begin{equation}\label{eq:energy}
	J(v;\Omega)=\int_\Omega\left(|\nabla v|^2+x_2\chi_{\{v>0\}}\right)\,\dx{x},
\end{equation}
where \(\chi_{\{v>0\}}\) denotes the characteristic function of
the set \(\{v>0\}\). We consider variational solutions of
\eqref{eq.problem} in the following sense.

\begin{definition}[Variational solution, cf. Definition 3.1 in \cite{VW11}]\label{def:variational}
A function \(u\in W^{1,2}_{\loc}(\Omega)\) is called a
\emph{variational solution} of \eqref{eq.problem} if \(u\in C^0(\Omega)\cap C^2(\Omega\cap\{u>0\})\), \(u\ge0\)
in \(\Omega\), \(u\equiv0\) in \(\Omega\cap\{x_2\le0\}\), and
\begin{equation}\label{eq:stationarity}
	\int_\Omega\left[
	(|\nabla u|^2+x_2\chi_{\{u>0\}})\diver\phi
	-2\,\nabla uD\phi\nabla u
	+\chi_{\{u>0\}}\phi_2\right]\dx{x}=0
\end{equation}
for every \(\phi=(\phi_1,\phi_2)\in C_c^1(\Omega;\R^2)\).
\end{definition}

\begin{remark}
In the domain variation identity \eqref{eq:stationarity}, the test functions are the compactly
supported functions
\(\phi=(\phi_1,\phi_2)\in C_c^1(\Omega;\R^2)\). Their supports may
intersect \(\Omega\cap\{x_2=0\}\), and no tangency condition is imposed
on them along this line. The formula \eqref{eq:stationarity} is nothing but the first variation of the energy functional \eqref{eq:energy}
with respect to domain variations: for every \(\phi\in C_c^1(\Omega;\R^2)\) and
every relatively compact open set \(\Omega'\Subset\Omega\) with
\(\operatorname{supp}\phi\subset\Omega'\),
\[
	0=-\frac{\dx{}}{\dx{\varepsilon}}J(u(x+\varepsilon\phi(x));\Omega')\Bigg|_{\varepsilon=0},
\]
independently of the choice of \(\Omega'\).
\end{remark}

Note that the \(C^2\) regularity in the positive phase, together with the domain
variation identity \eqref{eq:stationarity}, implies \(\Delta u=0\) in \(\{u>0\}\).
We also introduce the notion of weak solutions.

\begin{definition}[Weak solution]\label{def:weak}
A function \(u\) is called a \emph{weak solution} of
\eqref{eq.problem} if the following two conditions hold.
\begin{enumerate}[label=\textup{(\roman*)},leftmargin=*]
\item \(u\in W^{1,2}_{\loc}(\Omega)\cap C^0(\Omega)\),
\(u\ge0\) in \(\Omega\), \(u\equiv0\) in
\(\Omega\cap\{x_2\le0\}\), and \(u\) is harmonic in
\(\Omega\cap\{u>0\}\).

\item 
For every \(\tau>0\), the topological free boundary in
\(\Omega\cap\{x_2>\tau\}\) is locally the union of a relatively
open \(C^{2,\alpha}\) curve
\(\partial_{\mathrm{reg}}\{u>0\}\) and a
set of zero \(\mathcal H^1\) measure. Each point of this curve has an open
neighborhood \(U\) in which the positive phase lies on one side of
the curve, \(u\in C^1(U\cap\overline{\{u>0\}})\), and
\[
	|\nabla u|^2=x_2
	\qquad\text{on }U\cap\partial_{\mathrm{reg}}\{u>0\}.
\]
\end{enumerate}
\end{definition}

In Definition~\ref{def:weak}, the free boundary condition \textup{(ii)} 
is imposed above the line \(\{x_2=0\}\). A weak solution 
need not be variational. For example, let \(u(x)=a x_2^+\)
in \(B_1\), where \(a>0\). This function is harmonic in its
positive phase, and its free boundary above the line is empty,
so it satisfies Definition~\ref{def:weak}. However, for
\(\phi=(\phi_1,\phi_2)\in C_c^1(B_1;\R^2)\), the left hand side of
\eqref{eq:stationarity} equals
\[
	a^2\int_{B_1\cap\{x_2=0\}}\phi_2\,\dH{1},
\]
which need not vanish. Thus the weak formulation alone does
not imply stationarity across the line \(\{x_2=0\}\).

It should be noted that variational and weak solutions have different defining assumptions.
As shown above, a weak solution
need not satisfy the domain variation identity
\eqref{eq:stationarity}, hence need not be variational.
Definition~\ref{def:variational} imposes this identity without
assuming the free boundary regularity condition \textup{(ii)} in
Definition~\ref{def:weak}. By \cite[Lemma~3.4]{VW11}, a weak solution
satisfying \(|\nabla u|^2\le Cx_2^+\) locally in \(\Omega\) is a
variational solution.

In this paper the optimal gradient estimate of
Theorem~\ref{thm:growth} is proved for \emph{variational solutions},
and its core ingredient is the domain variation identity
\eqref{eq:stationarity}. We then obtain the blow-up alternatives
in Theorem~\ref{thm:stokes}\textup{(i)} for this class, and the
Stokes corner and free boundary regularity conclusions under the
additional hypotheses of parts~\textup{(ii)} and~\textup{(iii)}.

We next formally define stagnation points.
\begin{definition}[Stagnation points]
Let \(u\) be a variational solution of \eqref{eq.problem}. We call
\[
	\mathcal S^u=\left\lbrace x\in\Omega\cap\partial\{u>0\}:x_2=0
	\right\rbrace
\]
the set of stagnation points.
\end{definition}

Let \(u\) be a variational solution of \eqref{eq.problem}, and let
\(x^0\in\mathcal S^u\). The Bernoulli free boundary condition
\(|\nabla u|^2=x_2\) suggests the natural rescaling
\begin{equation}\label{eq.scaling}
	u_r(x):=\frac{u(x^0+rx)}{r^{3/2}}.
\end{equation}
Since \(x_2^0=0\), the domain variation identity \eqref{eq:stationarity} is preserved by this rescaling.
At points of the free boundary where
the Bernoulli condition \(|\nabla u|^2=x_2\) holds in the sense of a one sided \(C^1\) trace,
the condition rescales to \(|\nabla u_r|^2=x_2\) at the corresponding points of \(\partial\{u_r>0\}\).
At this stage, however, we lack the compactness needed to pass
to the limit as \(r\to0^+\) for the blow-up sequence defined in \eqref{eq.scaling}. 
The difficulty is the lack of growth and energy bounds that are
uniform in scale for \(r^{-3/2}u(x^0+rx)\). To overcome
this difficulty, the authors of \cite{VW11} imposed the so-called \textit{strong Bernstein estimate} as follows:
\begin{equation}\label{eq.strong-bernstein}
	|\nabla u|^2\le x_2^+\quad\text{ locally in }\Omega\cap\{u>0\},
\end{equation}
where \(x_2^+=\max\{x_2,0\}\). Under this gradient bound they proved the Stokes conjecture for weak solutions.
Throughout, the assumption that \(\{u=0\}\) has locally finitely many
connected components means that each stagnation point \(x^0\) has an
open neighborhood \(U\Subset\Omega\) such that
\(U\cap\{u=0\}\) has finitely many connected components.
\begin{theorem}[cf. Theorem A in \cite{VW11}]\label{thm.VW11}
Let \(u\) be a weak solution of \eqref{eq.problem} and suppose that \eqref{eq.strong-bernstein} holds. 
Suppose in addition that \(\{u=0\}\) has locally only finitely many connected 
components. If \(x^0=(x_1^0,0)\in\mathcal S^u\) is a stagnation point, then 
\(u_r\to\US\) as \(r\to0^+\) strongly in \(W^{1,2}_{\loc}(\R^2)\) and 
locally uniformly in \(\R^2\). Here \(\US\) is defined by
\begin{equation}\label{eq:stokes-profile}
	\US(\rho,\theta):=
	\begin{cases}
		\displaystyle\frac{\sqrt2}{3}\rho^{3/2}
		\cos\left[\frac32\left(\theta-\frac\pi2\right)\right],
		&\displaystyle\frac\pi6<\theta<\frac{5\pi}6,\\[3pt]
		0,&\text{otherwise},
	\end{cases}
\end{equation}
where \(x=\rho(\cos\theta,\sin\theta)\) and \(0\le\theta<2\pi\).
Moreover, in an open neighborhood of \(x^0\), the topological free boundary \(\partial\{u>0\}\) is the union of two
\(C^1\) graphs with right and left tangents at \(x^0\).
\end{theorem}

We make some comments regarding Theorem \ref{thm.VW11}.

\begin{remark}
Note that 
\(\{\US>0\}\) is an open sector with opening angle \(2\pi/3\), which 
confirms the conjecture proposed by G. Stokes (cf. Figure \ref{fig.stokes_corner}).
\begin{figure}[!htbp]
	\centering
	\begin{tikzpicture}[
		x=.93cm,
		y=.93cm,
		line cap=round,
		line join=round,
		every node/.style={font=\small},
		panel/.style={draw=black,line width=.35pt},
		level/.style={draw=black!45,line width=.35pt},
		freebdry/.style={draw=black,line width=1.05pt}
	]
		\def\xL{-2.15}
		\def\xR{ 2.15}
		\def\yB{-1.20}
		\def\yT{ 2.05}

		\begin{scope}[shift={(-2.65,2.80)}]
			\draw[panel]
				(-2.5,-1.8) rectangle (2.5,2.45);

			\draw[level] (-2.5,0) -- (2.5,0);

			\draw[freebdry]
				plot[domain=-2.5:0,samples=80]
				(\x,{-\x/sqrt(3)-.1025*\x*\x});

			\draw[freebdry]
				plot[domain=0:2.5,samples=80]
				(\x,{\x/sqrt(3)-.1025*\x*\x});

			\draw[dashed, line width=.5pt] (0,0) -- (30:2.3);
			\draw[dashed, line width=.5pt] (0,0) -- (150:2.3);

			\draw[line width=.35pt]
				(30:.26)
				arc[start angle=30,end angle=150,radius=.26];

			\node[
				font=\scriptsize,
				fill=white,
				inner sep=.4pt
			] at (.06,.47) {$120^\circ$};

			\node[font=\scriptsize, fill=white, inner sep=.4pt] at (30:2.05) {$\frac{\pi}{6}$};
			\node[font=\scriptsize, fill=white, inner sep=.4pt] at (150:2.05) {$\frac{5\pi}{6}$};

			\node at (0,1.55) {$\US>0$};
			\node at (0,-.92) {$\US=0$};

			\fill (0,0) circle (1.7pt);
		\end{scope}
	\end{tikzpicture}
	\caption{Stokes corner}
	\label{fig.stokes_corner}
\end{figure}
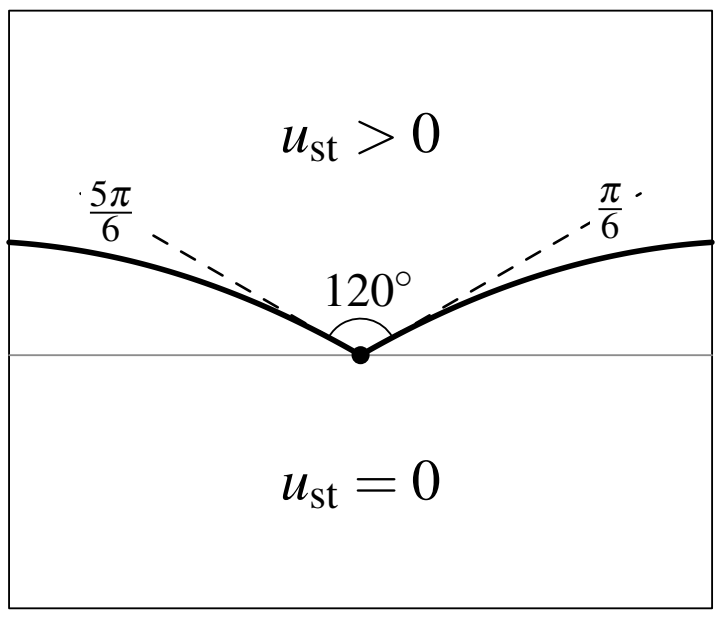
\end{remark}

\begin{remark}
The strong Bernstein estimate \eqref{eq.strong-bernstein} has two
roles in \cite{VW11}: it excludes cusp singularities
\cite[Lemma~4.4]{VW11} and enters the flatness argument used to
prove \(C^1\) regularity of the free boundary
\cite[Theorem~11.2]{VW11}.

The strong Bernstein
estimate \eqref{eq.strong-bernstein} has been verified in \cite{ST08} 
for waves of infinite depth.
In \cite{McC25}, McCurdy showed that it is not necessary for
excluding cusp singularities for weak 
solutions; the weaker growth bound
\eqref{eq:local-growth-intro} suffices for this purpose. More
recently, the cusp singularities were ruled out in all dimensions \(n\ge2\) under the growth
bound \(|\nabla u|^2\le Cx_n^+\), even for rotational flows
\cite{DPY26}.
\end{remark}

The main aim of this paper is to establish the local growth bound
\begin{equation}\label{eq:local-growth-intro}
	|\nabla u|^2\le Cx_2^+\quad\text{locally in }\Omega,
\end{equation}
for variational solutions of \eqref{eq.problem}. The principal role of
this bound is to supply the compactness required for the blow-up
analysis at stagnation points. In \cite{VW11}, a gradient bound of
this form is imposed as an a priori assumption on weak solutions; the
core contribution of the present paper is to prove that variational
solutions automatically satisfy such a bound. As an application, we
combine this estimate with the monotonicity, compactness, and frequency
theory of \cite{VW11} to establish Stokes corners and their one sided
\(C^1\) branches under the topological and interface hypotheses in
Theorem~\ref{thm:stokes}.

\section{Main results and ideas of the proof}

\subsection{Main results}

In this section, we state the main results of the paper. Throughout,
a constant \(C\) is called \emph{universal} if it is independent of
the solution, the center \(x^0\), and the radius \(R\). Our first result is an optimal local growth estimate for variational
solutions: \(|\nabla u|^2\) is controlled by the optimal factor \(x_2^+\) times the sum of the scaled local Dirichlet energy and an additive constant.

\begin{maintheorem}{I}[Optimal local growth]\label{thm:growth}
Let \(u\) be a variational solution of \eqref{eq.problem},
and let \(x^0=(x_1^0,0)\) with \(B_R(x^0)\Subset\Omega\)
for some \(R>0\). Then
\begin{equation}\label{eq:main-gradient}
	|\nabla u(x)|^2\le
	C\left(1+R^{-3}\int_{B_R(x^0)}|\nabla u|^2\dx{x}\right)x_2^+
	\quad\hbox{for a.e. }x\in B_{R/16}(x^0),
\end{equation}
where \(C>0\) is a universal constant.
\end{maintheorem}

Several comments are in order regarding Theorem~\ref{thm:growth}.

\begin{remark}
The optimality in Theorem~\ref{thm:growth} refers only to the exponents of
\(x_2^+\) in \eqref{eq:main-gradient}, not to the
multiplicative constant \(C\).
\end{remark}

\begin{remark}
Since \(x_2^+\) is bounded on \(B_{R/16}(x^0)\), \eqref{eq:main-gradient} yields
\(\nabla u\in L^\infty(B_{R/16}(x^0))\). Hence the continuous representative of \(u\) is Lipschitz on \(B_{R/16}(x^0)\). For this representative we have
\(u(x_1,0)=0\), since \(u=0\) in \(\{x_2<0\}\). For \(x_2>0\), integrating
\( |\partial_{x_2}u|\le|\nabla u| \) vertically from \((x_1,0)\) to
\((x_1,x_2)\) and using \eqref{eq:main-gradient}, we obtain the following
estimate for the solution itself after enlarging the universal constant \(C\),
\begin{equation}\label{eq:main-function}
	0\le u(x)\le
	C\left(1+R^{-3}\int_{B_R(x^0)}|\nabla u|^2\,\dx{x}\right)^{1/2}
	(x_2^+)^{3/2}
	\qquad\text{for all }x\in B_{R/16}(x^0).
\end{equation}
\end{remark}

\begin{remark}
The estimate \eqref{eq:main-function} also yields Fr\'echet
differentiability at every point
\(z=(z_1,0)\in B_{R/16}(x^0)\). Indeed, \(u(z)=0\), and for all sufficiently
small \(h=(h_1,h_2)\), we have
\[
	u(z+h)
	=O\!\left((h_2^+)^{3/2}\right)
	=O(|h|^{3/2})
	=o(|h|)
	\qquad\text{as }h\to0 .
\]
Hence \(u\) is Fr\'echet differentiable at \(z\), with \(Du(z)=0\).
\end{remark}

\begin{remark}
By Definition~\ref{def:variational},
\(u\in C^2(\Omega\cap\{u>0\})\). Hence \(\nabla u\) is continuous in the
positive phase \(\{u>0\}\), and the a.e. estimate
\eqref{eq:main-gradient} can be upgraded to the pointwise estimate
\begin{equation}\label{eq:main-gradient-pointwise}
	\begin{aligned}
		|\nabla u(x)|^2\le
		C\,\left(1+R^{-3}\int_{B_R(x^0)}|\nabla u|^2\,\dx{x}\right)x_2^+\\
		\qquad\qquad\qquad\qquad\text{for every }
		x\in B_{R/16}(x^0)\cap\{u>0\}.
	\end{aligned}
\end{equation}
Combined with the preceding function estimate and the resulting identity
\(Du=0\) on \(B_{R/16}(x^0)\cap\{x_2=0\}\), the pointwise bound
\eqref{eq:main-gradient-pointwise} also shows that
\[
	\nabla u(x)\longrightarrow0
	\qquad\text{as }\{u>0\}\ni x\to z
\]
for every \(z\in B_{R/16}(x^0)\cap\{x_2=0\}\). Thus \(\nabla u\) extends
continuously by zero from the positive phase to the line \(\{x_2=0\}\).
\end{remark}

Theorem~\ref{thm:growth} supplies the growth bound for the
geometric analysis of \cite{VW11}. 

\begin{remark}\label{rem:one-sided-trace}
Set \(\Gamma=\Omega\cap\partial\{u>0\}\cap\{x_2>0\}\).
We denote by \(\partial_{\mathrm{reg}}\{u>0\}\) the set of
points of \(\Gamma\) with a neighborhood \(U\) in which
\(\Gamma\) is an embedded \(C^{2,\alpha}\) curve,
\(U\cap\{u>0\}\) occupies one side of that curve, and
\(u\in C^1(U\cap\overline{\{u>0\}})\) satisfies
\(|\nabla u|^2=x_2\) on \(U\cap\Gamma\).
The points in
\(\Gamma\setminus\partial_{\mathrm{reg}}\{u>0\}\)
are said to be singular.
Here \(\mathcal H^1\) denotes one-dimensional Hausdorff measure.
With this definition of the regular part, the whole condition in
Definition~\ref{def:weak}\textup{(ii)} is equivalent to
\(\mathcal H^1(\Gamma\setminus
\partial_{\mathrm{reg}}\{u>0\})=0\).
\end{remark}

Above the line \(\{x_2=0\}\), the Bernoulli coefficient \(\sqrt{x_2}\) is positive
and analytic. In dimension two, the free boundary of a local
minimizer is regular there \cite[Theorem~8.3]{AC81}. For singular
points of the Bernoulli free boundary, see
\cite[Example~5.7]{AC81} and \cite[Section~3.1]{KW24}.

\Needspace{11\baselineskip}
Theorem~\ref{thm:growth} removes the need for an a priori gradient
bound. We now state the blow-up alternatives and the additional
hypotheses under which the Stokes corner and its two \(C^1\)
branches follow.

\begin{maintheorem}{II}[Stokes conjecture for variational solutions]
\label{thm:stokes}
Let \(u\) be a variational solution of \eqref{eq.problem}. Then the following hold.
\begin{enumerate}[label=\textup{(\roman*)},leftmargin=*]
\item At every stagnation point \(x^0\in\mathcal S^u\), the rescalings
\(u_r(x)=r^{-3/2}u(x^0+rx)\) satisfy exactly one of the following
two alternatives as \(r\to0^+\): either the Stokes corner asymptotic
\begin{equation}\label{eq:main-blowup}
	u_r\longrightarrow\US,\qquad
	\chi_{\{u_r>0\}}\longrightarrow\chi_{\{\US>0\}},
\end{equation}
or the full-density asymptotic
\begin{equation}\label{eq:main-flat-blowup}
	u_r\longrightarrow0,\qquad
	\chi_{\{u_r>0\}}\longrightarrow\chi_{\{x_2>0\}}.
\end{equation}
In both cases \(u_r\to u_0\) strongly in \(W^{1,2}_{\loc}(\R^2)\)
and locally uniformly, and \(\chi_{\{u_r>0\}}\) converges in
\(L^1_{\loc}(\R^2)\). These convergences hold for the whole family
as \(r\to0^+\). Here \(u_0\) denotes the corresponding blow-up
limit, and \(\US\) is defined in \eqref{eq:stokes-profile}.

\item Suppose that \(\{u=0\}\) has locally finitely many connected
components and that
\(\mathcal H^1(\Gamma\setminus\partial_{\mathrm{reg}}\{u>0\})=0\).
Then \(\mathcal S^u\) is locally finite and \eqref{eq:main-blowup}
holds at every stagnation point. Here \(\Gamma:=\Omega\cap\partial\{u>0\}\cap\{x_2>0\}\)
and \(\partial_{\mathrm{reg}}\{u>0\}\) denotes the regular part of the free boundary
defined in Remark~\ref{rem:one-sided-trace}.

\item Suppose that the hypotheses of \textup{(ii)} hold and that
\(\Gamma\setminus\partial_{\mathrm{reg}}\{u>0\}=\varnothing\).
Then near each \(x^0\in\mathcal S^u\) the topological free boundary
consists of two graph branches meeting at \(x^0\). Each is \(C^1\)
up to \(x^0\), with slope \(-1/\sqrt3\) on the left and \(1/\sqrt3\)
on the right.
\end{enumerate}
\end{maintheorem}

\begin{remark}
The proof of Theorem~\ref{thm:stokes} combines our gradient and interface
estimates with the monotonicity, compactness, density, and frequency theory
of \cite{VW11}. The compactness theory in
\cite[Sections~3--8]{VW11} is formulated for the variational class in
Definition~\ref{def:variational}.

\begin{itemize}
\item The essential difference between Theorem~\ref{thm:stokes} and Theorem~A of \cite{VW11} concerns the solution class. In \cite{VW11}, solutions are assumed to be weak solutions satisfying the strong Bernstein estimate \eqref{eq.strong-bernstein}. By \cite[Lemma~3.4]{VW11}, such solutions are variational, and \cite[Sections~3--8]{VW11} develop the compactness and limit theory for variational solutions under the growth condition. Here we work directly with variational solutions and prove the optimal growth bound of Theorem~\ref{thm:growth}, so that our solutions automatically satisfy the compactness conditions of \cite[Sections~3--8]{VW11}; this yields the convergence in part~\textup{(i)}.

The exclusion of cusp singularities, namely the case \(u_r\to0\) and \(\chi_{\{u_r>0\}}\to0\), is proved in \cite{VW11} for weak solutions under the strong Bernstein estimate; Proposition~\ref{prop:no-cusps-variational} proves it for every variational solution.

\item The exclusion of full-density singularities is the same as in \cite{VW11}; here it is stated as Proposition~\ref{prop:no-full-density}.

\item In \cite[Theorem~11.2]{VW11}, the bound
\(|\nabla u|^2\le x_2^+\) is used in the flatness argument after
the Stokes profile has been identified. Here Theorem~\ref{thm:growth}
gives a bound with a general constant \(C\). For the flatness
argument in part~\textup{(iii)}, we recover the coefficient one
asymptotically (cf. Lemma~\ref{lem:pressure-excess}): for each fixed \(L>0\),
\[
	|\nabla u_r|^2\le(1+\varepsilon_r)x_2
	\quad\text{in }B_L\cap\{u_r>0\},\qquad
	\varepsilon_r=\varepsilon(r;L)\to0.
\]
\end{itemize}
\end{remark}

A central inspiration for our approach is the recent elegant work of
De~Philippis, Hirsch, and Nahon \cite{DHN26}. In the context of
stationary Bernoulli free boundary problems, they introduced a
Hopf differential argument, combined with the Beurling transform
and their abstract approximation theorem, to establish optimal
Lipschitz regularity.
Our strategy is strongly inspired by this work. Nevertheless, our
starting point is the gravity water wave problem, whose
domain variation identity contains the weighted Bernoulli term
\(x_2\chi_{\{u>0\}}\). Consequently, the unweighted Lipschitz
result cannot be invoked directly.

\subsection{Ideas of the proof}

We write \(B_s(\mathbf{b})\) for the disk of radius \(s\) centered at \(\mathbf{b}=(b_1,b_2)\)
and \(B_1=B_1(0)\). For a measurable set \(E\), \(|E|\) denotes
its area. We write \(I\) for the \(2\times2\) identity and
\(U\Subset\Omega\) when \(\overline U\) is compact in \(\Omega\).
In complex coordinates, we write
\[
	z=x_1+ix_2,\qquad
	\partial_z=\tfrac12(\partial_{x_1}-i\partial_{x_2}),\qquad
	\partial_{\bar z}=\tfrac12(\partial_{x_1}+i\partial_{x_2}).
\]

The domain variation identity \eqref{eq:stationarity} gives the
distributional equation
\[
	\partial_{\bar z}\Phi=-x_2\partial_z\chi_{\{u=0\}},
	\qquad \Phi=(u_{x_1}-iu_{x_2})^2.
\]
The Dirichlet energy controls the \(L^1\) norm of \(\Phi\) because
\(|\Phi|=|\nabla u|^2\). To obtain \(L^2\) control, we use the
Cauchy and Beurling transforms to write \(\Phi=h_B-E_B\) on a disk
\(B\), where \(h_B\) is holomorphic and \(E_B\in L^2(B)\).
Applying the mean-value inequality for subharmonic functions to
\(|h_B|\) yields a bound for \(h_B\) on a smaller disk by its \(L^1\) norm. Subtracting \(E_B\) gives the local \(L^2\) bound
for \(\Phi\) in Lemma~\ref{lem:approximation}.

On a disk \(B_r(\mathbf{a})\) centered on the line \(\{x_2=0\}\), the
\(L^2\) approximation error is at most \(Cr^2\) because
\(|x_2|\le r\) there. Since \(\Phi=0\) below the line, the same
error bounds the norm of
\(h_{B_r(\mathbf{a})}\) on the lower half disk. The holomorphic estimate
\eqref{eq:holomorphic-observation} then gives
\[
	\|\Phi\|_{L^2(B_{r/2}(\mathbf{a}))}
	\le\varepsilon\|\Phi\|_{L^2(B_r(\mathbf{a}))}+C_\varepsilon r^2.
\]
For \(\varepsilon<1/4\), dyadic iteration yields
\(\|\Phi\|_{L^2(B_r(\mathbf{a}))}\le Cr^2\), with \(C\) controlled by
the local Dirichlet energy.

Let us address the main differences between our formulation and
\cite{DHN26}. The first appears in the equation satisfied by the Hopf
differential \(\Phi\). For the unweighted Bernoulli problem, writing
\(\Phi_{\mathrm{unw}}=4(\dz u)^2\), the stationarity condition gives
\[
	\dbz\Phi_{\mathrm{unw}}
	=-\dz\chi_{\{u=0\}};
\]
see \cite[Lemma~18]{DHN26}. In our setting, the corresponding equation is instead (cf. Lemma~\ref{lem:hopf})
\[
	\dbz\Phi=-x_2\dz\chi_{\{u=0\}}.
\]
Thus, the constant coefficient defect in the unweighted problem is
replaced by a defect involving the variable weight \(x_2\).
The second difference concerns the correction. For the unweighted
problem, a Beurling correction alone suffices, since
\(\dbz\BT\chi_{\{u=0\}}=\dz\chi_{\{u=0\}}\). In our setting, the
analogous weighted Beurling term satisfies
\[
	\dbz\BT\bigl(x_2\chi_{\{u=0\}}\bigr)
	=\dz\bigl(x_2\chi_{\{u=0\}}\bigr)
	=x_2\dz\chi_{\{u=0\}}
	+(\dz x_2)\chi_{\{u=0\}}.
\]
The last term is an additional defect produced by the product
rule, and it necessitates an additional Cauchy transform
correction. We therefore introduce the correction
\[
	E_B=\BT\bigl(x_2\chi_{\{u=0\}\cap B}\bigr)
	-\CT\bigl((\dz x_2)\chi_{\{u=0\}\cap B}\bigr),
	\qquad h_B=\Phi+E_B.
\]
Since \(\dbz\CT g=g\), the Cauchy term exactly cancels
\((\dz x_2)\chi_{\{u=0\}}\), and hence \(h_B\) is holomorphic
(Lemma~\ref{lem:holomorphic-correction}). This additional Cauchy
correction is the new ingredient needed to adapt the
Hopf--Beurling approach to the gravity water wave problem.

Once \(h_B\) is holomorphic, the second hypothesis of the abstract
approximation theorem follows from the standard interior estimate
for holomorphic functions, with a universal constant. Moreover,
on every subdisk of \(B_r(\mathbf{a})\), the squared approximation error is
bounded by \(Cr^2\) times the measure of the zero set of \(\Phi\)
in that subdisk. The zero phase below the line \(\{x_2=0\}\) also
yields the essential decay
\(\|\Phi\|_{L^2(B_r(\mathbf{a}))}\le Cr^2\). The scaled approximation
theorem (Theorem~\ref{thm:DHN}) therefore gives
\[
	\|\Phi\|_{L^\infty(B_{r/2}(\mathbf{a}))}
	\le C\bigl(r+r^{-1}\|\Phi\|_{L^2(B_r(\mathbf{a}))}\bigr)
	\le Cr.
\]
Choosing \(\mathbf{a}=(x_1,0)\) and \(r=4x_2\) gives
Theorem~\ref{thm:growth} after rescaling.

For Theorem~\ref{thm:stokes}\textup{(i)}, the gradient bound and
weighted perimeter estimate give the compactness needed for
the Weiss analysis of \cite{VW11}. The weighted isoperimetric
argument of \cite{DPY26}, given in
Proposition~\ref{prop:no-cusps-variational}, excludes zero density
in the variational class. In part~\textup{(ii)}, the normalized frequency profile has an
interior nodal ray. The finite component condition and phase
convergence give, by radial projection, a positive lower bound
for the interface length on a fixed disk above the line
\(\{x_2=0\}\). The vanishing of the \(\mathcal H^1\) measure of the singular
set and the one sided Bernoulli trace turn this into a positive
lower bound for \(\Delta u_r\) on that disk, contradicting
\(u_r\to0\).

For part~\textup{(iii)}, extending \((|\nabla u_r|^2-x_2)^+\)
by zero gives a subharmonic function. Comparison with a harmonic
barrier then yields the asymptotically sharp gradient bound. We localize the interface near the two Stokes
rays on annuli and apply De~Silva's flatness theorem
\cite{DS11} to control the normals. Counting crossings of vertical
segments gives the two graph branches.

\section{The Hopf differential formulation}

After translating a point \(x^0=(x_1^0,0)\) of the line
\(\{x_2=0\}\) to the origin and rescaling by
\(u_r(x)=r^{-3/2}u(x^0+rx)\), we work in \(B_1=B_1(0)\) and
denote the rescaled function by \(u\).

We use the complex notation introduced above. The Hopf differential
of \(u\) is defined by 
\begin{equation}\label{eq:hopf-definition}
	\Phi=(u_{x_1}-iu_{x_2})^2.
\end{equation}
Its modulus satisfies \( |\Phi|=|\nabla u|^2 \). In particular,
\(\Phi\in L^1(B_1)\) whenever \(u\in W^{1,2}(B_1)\). Moreover, \(\Phi=0\) almost everywhere on \(\{u=0\}\), since
\(\nabla u=0\) almost everywhere on \(\{u=0\}\).

\begin{lemma}\label{lem:hopf}
Let \(u\in W^{1,2}(B_1)\) be nonnegative. Assume that \(u\) satisfies the domain variation identity
\eqref{eq:stationarity} on \(\Omega=B_1\) for all
\(\phi\in C_c^1(B_1;\R^2)\). Then the Hopf differential \eqref{eq:hopf-definition} satisfies
\begin{equation}\label{eq:hopf}
	\dbz\Phi=-x_2\dz\chi_{\{u=0\}}\qquad\hbox{ in }B_1,
\end{equation}
in the sense of distributions.
\end{lemma}

\begin{proof}
Set \(A=u_{x_1}^2-u_{x_2}^2\), \(B=2u_{x_1}u_{x_2}\), and \(q=x_2\chi_{\{u>0\}}\).
Then \(A,B,q\in L^1(B_1)\) and \(\Phi=A-iB\). The matrix
appearing in \eqref{eq:stationarity} is
\[
	(|\nabla u|^2+q)I-2\nabla u\otimes\nabla u
	=\begin{pmatrix}-A+q&-B\\-B&A+q\end{pmatrix},
\]
where \((\nabla u\otimes\nabla u)_{ij}=u_{x_i}u_{x_j}\).
Taking \(\phi=(\eta,0)\) and \((0,\eta)\), respectively, with
\(\eta\in C_c^1(B_1)\), in \eqref{eq:stationarity} gives
\begin{equation}\label{eq:stationarity-components}
	A_{x_1}+B_{x_2}=q_{x_1},\qquad
	A_{x_2}-B_{x_1}=\chi_{\{u>0\}}-q_{x_2}
\end{equation}
in distributions. Hence
\begin{align*}
	\dbz\Phi
	&=\frac12(A_{x_1}+B_{x_2})
	+\frac i2(A_{x_2}-B_{x_1})\\
	&=\dz q+\frac i2\chi_{\{u>0\}}.
\end{align*}
Since \(\dz x_2=-i/2\), the distributional product rule yields
\[
	\dz q=-\frac i2\chi_{\{u>0\}}
	+x_2\dz\chi_{\{u>0\}}.
\]
Thus \(\dbz\Phi=x_2\dz\chi_{\{u>0\}}\), and \eqref{eq:hopf} follows from
\(\chi_{\{u>0\}}=1-\chi_{\{u=0\}}\) a.e.
\end{proof}

\begin{remark}
The second identity in \eqref{eq:stationarity-components} uses
\(\phi=(0,\eta)\) in \eqref{eq:stationarity}, with arbitrary
\(\eta\in C_c^1(B_1)\). If the test functions were required to be
tangent to \(\{x_2=0\}\), then \(\eta\) would vanish there, leaving
a possible distribution supported on the line:
\[
	A_{x_2}-B_{x_1}
	=\chi_{\{u>0\}}-q_{x_2}+\sigma.
\]
Here \(\sigma\) is a distribution of the form
\(\sigma=T(x_1)\delta_{x_2=0}\) for some distribution \(T\) on
\(B_1\cap\{x_2=0\}\), and \(\sigma\) has no prescribed
sign.
The Hopf equation would then contain \(\frac{i}{2}\sigma\).
The approximation and decay arguments below require
\eqref{eq:hopf} on disks crossing the line.
\end{remark}

\section{Holomorphic approximation and
\texorpdfstring{\(L^2\)}{L2} bounds for the Hopf differential}

\enlargethispage{2pt}

We construct holomorphic approximants of \(\Phi\) whose \(L^2\)
errors are controlled by the measure of the zero phase. We also
establish the local \(L^2\) integrability needed to apply the
approximation theorem of De~Philippis, Hirsch, and Nahon
\cite[Theorem~1]{DHN26}.

On each disk \(B=B_s(\mathbf{b})\subset B_1\), the Cauchy and Beurling
transforms give a correction \(E_B\in L^2(B)\) such that
\(h_B=\Phi+E_B\) is holomorphic. Their bounds control the
approximation error by
\((s+\sup_B|x_2|)|B\cap\{u=0\}|^{1/2}\).
The factor involving \(x_2\) will give the decay near the line \(\{x_2=0\}\).

The same construction improves the integrability of \(\Phi\).
Applying the mean-value inequality for subharmonic functions to
\(|h_B|\) yields a bound for \(h_B\) on smaller disks by its \(L^1\) norm. Since \(\Phi=h_B-E_B\) and \(E_B\in L^2(B)\),
we obtain local \(L^2\) control of \(\Phi\).

\begin{lemma}\label{lem:approximation}
Let \(u\in W^{1,2}(B_1)\) be nonnegative and satisfy the domain
variation identity \eqref{eq:stationarity} on \(\Omega=B_1\) for every
\(\phi\in C_c^1(B_1;\R^2)\). Then, for every disk
\(B=B_s(\mathbf{b})\subset B_1\), there exists a holomorphic function \(h_B\)
on \(B\) such that 
\begin{equation}\label{eq:approximation}
	\|\Phi-h_B\|_{L^2(B)}
	\le C\bigl(s+\sup_B|x_2|\bigr)
	|B\cap\{u=0\}|^{1/2},
\end{equation}
where \(C>0\) is universal. Moreover,
\[
	\Phi\in L^2(B_{3/4}),\qquad
	\|\Phi\|_{L^2(B_{3/4})}
	\le C\left(1+\int_{B_1}|\nabla u|^2\dx{x}\right).
\]
In particular, \(h_B\in L^2(B)\) for every disk \(B\subset B_{3/4}\).
\end{lemma}

Lemma~\ref{lem:approximation} follows from
Lemmas~\ref{lem:holomorphic-correction} and~\ref{lem:integrability-height}
below. We first recall the two transforms used in their proofs.
For a compactly supported function
\(g\in L^2(\R^2)\), define the Cauchy and Beurling transforms (see \cite[Definition 19]{DHN26} for the Beurling transform) by
\begin{equation}\label{eq:transforms}
	\CT g(z)=\frac1\pi\int_{\R^2}\frac{g(\zeta)}{z-\zeta}
	\,\dx{\zeta},
	\qquad
	\BT g=\dz\CT g.
\end{equation}
The convolution \(\CT g\) in \eqref{eq:transforms} is locally
integrable, and \(\BT g\) is
defined in distributions. Since \(1/(\pi z)\) is a fundamental
solution of \(\dbz\),
\begin{equation}\label{eq:transform-identities}
	\dbz\CT g=g,
	\qquad
	\dbz\BT g=\dz g.
\end{equation}
Thus, if \(\Omega\subset\mathbb C\) is open,
\(f\in L^1_{\loc}(\Omega)\), \(g\in L^2(\R^2)\) is compactly
supported, and \(\dbz f=g\) in the sense of distributions in
\(\Omega\), then \(f-\CT g\) has a holomorphic representative
in \(\Omega\).
The Beurling transform is an \(L^2\) isometry: its Fourier
symbol \footnote{ We use the convention
\(\widehat f(\xi)=\int_{\R^2}e^{-ix\cdot\xi}f(x)\,\dx{x}\).
For an operator \(T\) acting by multiplication on the Fourier side,
its symbol \(m_T\) is defined by
\(\widehat{Tf}=m_T\widehat f\). Thus \(\dz\) and \(\dbz\) have symbols
\(\frac{i}{2}(\xi_1-i\xi_2)\) and
\(\frac{i}{2}(\xi_1+i\xi_2)\), respectively.} is
\[
	\frac{\xi_1-i\xi_2}{\xi_1+i\xi_2},
	\qquad \xi\ne0,
\]
which has modulus one. Plancherel's theorem therefore gives
\begin{equation}\label{eq:beurling-isometry}
	\|\BT g\|_{L^2(\R^2)}=\|g\|_{L^2(\R^2)}.
\end{equation}

For the Cauchy transform, suppose that
\(\supp g\subset\overline{B_s(\mathbf{b})}\). Since \(|z-\zeta|\le2s\)
for \(z\in B_s(\mathbf{b})\) and \(\zeta\in\supp g\), we obtain
\[
	|\CT g(z)|
	\le\int_{\R^2}
	\frac{\chi_{B_{2s}}(z-\zeta)}{\pi|z-\zeta|}
	|g(\zeta)|\,\dx{\zeta}.
\]
The truncated kernel has integral
\[
	\int_{\R^2}\frac{\chi_{B_{2s}}(w)}{\pi|w|}\,\dx{w}
	=\int_{B_{2s}}\frac1{\pi|w|}\,\dx{w}
	=\frac1\pi\int_0^{2s}\int_0^{2\pi}\,\dx{\theta}\,\dx{\rho}
	=4s.
\]
Young's inequality gives
\begin{equation}\label{eq:cauchy-local-bound}
	\|\CT g\|_{L^2(B_s(\mathbf{b}))}
	\le4s\|g\|_{L^2(\R^2)}.
\end{equation}

For the unit disk, the transforms can be computed explicitly.
\begin{example}
For \(g=\chi_{B_1}\), direct computation gives
\[
	\CT g(z)=
	\begin{cases}
		\bar z,& |z|<1,\\[2pt]
		z^{-1},& |z|>1,
	\end{cases}
	\qquad
	\BT g(z)=
	\begin{cases}
		0,& |z|<1,\\[2pt]
		-z^{-2},& |z|>1.
	\end{cases}
\]
In particular, \(\dbz\CT g=g\), and
\(\|\BT g\|_{L^2(\R^2)}^2=\pi=\|g\|_{L^2(\R^2)}^2\).
\end{example}

\begin{lemma}[Holomorphic correction]\label{lem:holomorphic-correction}
Let \(u\) satisfy the hypotheses of Lemma~\ref{lem:approximation}.
For every disk \(B=B_s(\mathbf{b})\subset B_1\), define
\begin{equation}\label{eq:holomorphic-correction}
	\begin{split}
		E_B&=\BT\bigl(x_2\chi_{\{u=0\}\cap B}\bigr)
		-\CT\bigl((\dz x_2)\chi_{\{u=0\}\cap B}\bigr),\\
		h_B&=\Phi+E_B,
	\end{split}
\end{equation}
where \(x_2\chi_{\{u=0\}\cap B}\) and \((\dz x_2)\chi_{\{u=0\}\cap B}\) 
are extended by zero outside \(B\). Then \(h_B\) has
a holomorphic representative in \(B\), and
\[
	\|\Phi-h_B\|_{L^2(B)}
	=\|E_B\|_{L^2(B)}
	\le\bigl(\sup_B|x_2|+2s\bigr)
	|B\cap\{u=0\}|^{1/2}.
\]
In particular, we obtain
\[
	\|\Phi-h_B\|_{L^2(B)}
	\le2\bigl(s+\sup_B|x_2|\bigr)
	|B\cap\{u=0\}|^{1/2}.
\]
Moreover, if \(\Phi\in L^2(B)\), then \(h_B\in L^2(B)\).
\end{lemma}

\begin{proof}
The arguments of \(\BT\) and \(\CT\) in
\eqref{eq:holomorphic-correction} are compactly supported in
\(\overline B\) and belong to \(L^2(\R^2)\), since
\(x_2\) is bounded on \(B\) and
\[
	\dz x_2
	=\tfrac12(\partial_{x_1}x_2-i\partial_{x_2}x_2)
	=-\frac i2.
\]
By \eqref{eq:beurling-isometry}, we estimate
\begin{align*}
	\left\|\BT\bigl(x_2\chi_{\{u=0\}\cap B}\bigr)\right\|_{L^2(B)}
	&\le\left\|\BT\bigl(x_2\chi_{\{u=0\}\cap B}\bigr)
	\right\|_{L^2(\R^2)}\\
	&=\left(\int_{B\cap\{u=0\}}x_2^2\,\dx{x}\right)^{1/2}\\
	&\le\sup_B|x_2|\,|B\cap\{u=0\}|^{1/2}.
\end{align*}
Likewise, \eqref{eq:cauchy-local-bound} gives
\begin{align*}
	\left\|\CT\bigl((\dz x_2)\chi_{\{u=0\}\cap B}\bigr)
	\right\|_{L^2(B)}
	&\le4s\left(\int_{B\cap\{u=0\}}|\dz x_2|^2\,\dx{x}\right)^{1/2}\\
	&=2s\,|B\cap\{u=0\}|^{1/2}.
\end{align*}
Adding the two bounds, we obtain
\begin{equation}\label{eq:correction-bound}
	\|E_B\|_{L^2(B)}
	\le\bigl(\sup_B|x_2|+2s\bigr)
	|B\cap\{u=0\}|^{1/2}.
\end{equation}
For \(\eta\in C_c^\infty(B)\), \eqref{eq:transform-identities} gives
\begin{align*}
	\langle\dbz E_B,\eta\rangle
	&=-\int_{B\cap\{u=0\}}x_2\,\dz\eta\,\dx{x}
	-\int_{B\cap\{u=0\}}(\dz x_2)\eta\,\dx{x}\\
	&=-\int_B\chi_{\{u=0\}}\,\dz(x_2\eta)\,\dx{x}\\
	&=\langle x_2\dz\chi_{\{u=0\}},\eta\rangle .
\end{align*}
Here we used \(\supp\eta\subset B\) and the product rule. Thus
\[
	\dbz E_B=x_2\dz\chi_{\{u=0\}}
\]
in distributions on \(B\). By \eqref{eq:hopf}, we compute
\[
	\dbz h_B
	=\dbz\Phi+\dbz E_B
	=-x_2\dz\chi_{\{u=0\}}+x_2\dz\chi_{\{u=0\}}
	=0
	\qquad\hbox{in }B.
\]
Moreover, \(\Phi\in L^1(B)\), so H\"older's inequality and
\eqref{eq:correction-bound} give
\[
	\|h_B\|_{L^1(B)}
	\le\|\Phi\|_{L^1(B)}
	+|B|^{1/2}\|E_B\|_{L^2(B)}
	<\infty.
\]
Weyl's lemma (see \cite[Corollary~2.2.1]{Jost13}) applied to \(\Delta h_B=4\dz\dbz h_B=0\) gives a
smooth representative. Since \(\dbz h_B=0\), it is holomorphic.

The error bounds follow from \(\Phi-h_B=-E_B\),
\eqref{eq:correction-bound}, and
\[
	\sup_B|x_2|+2s\le2\bigl(s+\sup_B|x_2|\bigr).
\]
If \(\Phi\in L^2(B)\), then
\[
	\|h_B\|_{L^2(B)}
	\le\|\Phi\|_{L^2(B)}+\|E_B\|_{L^2(B)}
	<\infty,
\]
which proves the last assertion.
\end{proof}

The correction estimate and the mean-value inequality for subharmonic
functions give local \(L^2\) integrability of \(\Phi\).

\begin{lemma}[Integrability and estimates at the height scale]
\label{lem:integrability-height}
Let \(u\) satisfy the hypotheses of Lemma~\ref{lem:approximation}, and
let \(h_B\) be the approximant constructed in
Lemma~\ref{lem:holomorphic-correction}. Then
\[
	\Phi\in L^2(B_{3/4}),\qquad
	\|\Phi\|_{L^2(B_{3/4})}
	\le C\left(1+\int_{B_1}|\nabla u|^2\,\dx{x}\right),
\]
where \(C>0\) is universal. In particular, \(h_B\) belongs to \(L^2(B)\) for every
disk \(B\subset B_{3/4}\).
Furthermore, if \(\mathbf{a}=(a_1,0)\) and
\(B_s(\mathbf{b})\subset B_r(\mathbf{a})\subset B_{3/4}\), then
\[
	\int_{B_s(\mathbf{b})}|\Phi-h_{B_s(\mathbf{b})}|^2\,\dx{x}
	\le9r^2|B_s(\mathbf{b})\cap\{\Phi=0\}|.
\]
\end{lemma}

\begin{proof}
Apply \eqref{eq:correction-bound} with \(B=B_1\). Since
\(\sup_{B_1}|x_2|=1\) and \(|B_1\cap\{u=0\}|\le\pi\), we obtain
\[
	\|E_{B_1}\|_{L^2(B_1)}
	\le3|B_1\cap\{u=0\}|^{1/2}
	\le3\sqrt\pi.
\]
Using \(|\Phi|=|\nabla u|^2\) and H\"older's inequality, we obtain
\begin{align*}
	\|h_{B_1}\|_{L^1(B_1)}
	&\le\|\Phi\|_{L^1(B_1)}+\|E_{B_1}\|_{L^1(B_1)}\\
	&\le\int_{B_1}|\nabla u|^2\,\dx{x}
	+|B_1|^{1/2}\|E_{B_1}\|_{L^2(B_1)}\\
	&\le\int_{B_1}|\nabla u|^2\,\dx{x}+3\pi.
\end{align*}

For \(z\in B_{3/4}\), the disk \(B_{1/8}(z)\) lies in \(B_1\).
Applying the mean-value inequality for subharmonic functions to
\(|h_{B_1}|\) gives
\[
	\begin{aligned}
		|h_{B_1}(z)|
		&\le\frac1{|B_{1/8}|}
		\int_{B_{1/8}(z)}|h_{B_1}(\zeta)|\,\dx{\zeta}\\
		&\le\frac{64}{\pi}\|h_{B_1}\|_{L^1(B_1)}\\
		&\le C\left(1+\int_{B_1}|\nabla u|^2\,\dx{x}\right).
	\end{aligned}
\]
Since \(\Phi=h_{B_1}-E_{B_1}\), it follows that
\begin{equation}\label{eq:initial-hopf-bound}
	\begin{aligned}
		\|\Phi\|_{L^2(B_{3/4})}
		&\le\|h_{B_1}\|_{L^2(B_{3/4})}
		+\|E_{B_1}\|_{L^2(B_{3/4})}\\
		&\le|B_{3/4}|^{1/2}
		\|h_{B_1}\|_{L^\infty(B_{3/4})}
		+\|E_{B_1}\|_{L^2(B_1)}\\
		&\le C\left(1+\int_{B_1}|\nabla u|^2\,\dx{x}\right).
	\end{aligned}
\end{equation}
In particular, \(h_B\in L^2(B)\) for \(B\subset B_{3/4}\), by
Lemma~\ref{lem:holomorphic-correction}.

Finally, let \(\mathbf{a}=(a_1,0)\) and
\(B_s(\mathbf{b})\subset B_r(\mathbf{a})\subset B_{3/4}\). Then \(s\le r\) and
\(\sup_{B_s(\mathbf{b})}|x_2|\le r\). Squaring \eqref{eq:correction-bound} yields
\begin{equation}\label{eq:approximation-height}
	\begin{aligned}
		\int_{B_s(\mathbf{b})}|\Phi-h_{B_s(\mathbf{b})}|^2\,\dx{x}
		&\le\bigl(\sup_{B_s(\mathbf{b})}|x_2|+2s\bigr)^2
		|B_s(\mathbf{b})\cap\{u=0\}|\\
		&\le9r^2|B_s(\mathbf{b})\cap\{u=0\}|\\
		&\le9r^2|B_s(\mathbf{b})\cap\{\Phi=0\}|.
	\end{aligned}
\end{equation}
The last inequality follows from \(\Phi=0\) almost everywhere on
\(\{u=0\}\).
\end{proof}

\begin{remark}
The estimate \eqref{eq:initial-hopf-bound} also gives
\[
	\|\nabla u\|_{L^4(B_{3/4})}^2
	=\|\Phi\|_{L^2(B_{3/4})}<\infty.
\]
Thus \(\nabla u\in L^4(B_{3/4})\), and hence
\(u\in W^{1,4}(B_{3/4})\). It then follows from Morrey's inequality
\cite[Theorem~11.1.3]{Jost13} that \(u\) is H\"older continuous
with exponent \(1/2\), and in particular has a continuous
representative on this disk.
\end{remark}

Combining Lemmas~\ref{lem:holomorphic-correction}
and~\ref{lem:integrability-height} proves Lemma~\ref{lem:approximation}.

\begin{remark}
In the approximation theorem of De~Philippis, Hirsch, and Nahon
\cite[Theorem~1, estimate~(1.7)]{DHN26}, the approximation error is
controlled by the measure of the zero set. In the scaled form stated in
Theorem~\ref{thm:DHN}, this hypothesis requires
\[
	\int_{B_s(\mathbf{b})}|\Phi-h_{B_s(\mathbf{b})}|^2\,\dx{x}
	\le\Lambda^2|B_s(\mathbf{b})\cap\{\Phi=0\}|
\]
for every \(B_s(\mathbf{b})\subset B_r(\mathbf{a})\subset B_{3/4}\), where \(a_2=0\).
Lemma~\ref{lem:integrability-height} verifies precisely this condition
(see \eqref{eq:approximation-height}) with \(\Lambda=3r\). It also supplies the required \(L^2\) integrability
of \(\Phi\) and of the approximants \(h_{B_s(\mathbf{b})}\).
\end{remark}

\section{Proof of gradient estimates}
\label{sec:decay}

We first estimate the Hopf differential on disks centered on the line \(\{x_2=0\}\), then apply \cite[Theorem~1]{DHN26} to obtain a
supremum bound. Choosing disks with radius proportional to the height yields
Theorem~\ref{thm:growth}.

\begin{lemma}\label{lem:decay}
Let \(u\in W^{1,2}(B_1)\) be nonnegative and vanish almost everywhere
on \(B_1\cap\{x_2\le0\}\), and suppose that \(u\) satisfies
\eqref{eq:stationarity} on \(\Omega=B_1\) for every
\(\phi\in C_c^1(B_1;\R^2)\). Set
\(\Phi=(u_{x_1}-iu_{x_2})^2\). Then, for
\(\mathbf{a}=(a_1,0)\), \(|a_1|\le1/4\), and \(0<r\le1/4\),
\begin{equation}\label{eq:decay}
	\|\Phi\|_{L^2(B_r(\mathbf{a}))}
	\le C\left(1+\int_{B_1}|\nabla u|^2\dx{x}\right)r^2.
\end{equation}
Here \(C>0\) is a universal constant.
\end{lemma}

\begin{proof}
We first record the interior estimate
\begin{equation}\label{eq:holomorphic-interior}
	s\|\nabla h\|_{L^\infty(B_{s/2}(\mathbf{b}))}
	\le C\left(\barint_{B_s(\mathbf{b})}
	|h-(h)_{B_s(\mathbf{b})}|^2\dx{x}\right)^{1/2}
\end{equation}
for every holomorphic \(h\in L^2(B_s(\mathbf{b}))\).
Indeed, applying the interior gradient estimate for harmonic functions
\cite[Corollary~2.2.7]{Jost13} to the real and imaginary parts of
\(h-(h)_{B_s(\mathbf{b})}\) on balls centered at points of
\(B_{s/2}(\mathbf{b})\) and contained in \(B_s(\mathbf{b})\), and using
the mean value inequality for subharmonic functions and H\"older's
inequality, gives
\[
	s|\nabla h(z)|
	\le C s^{-2}\|h-(h)_{B_s(\mathbf{b})}\|_{L^1(B_s(\mathbf{b}))}
	\le C\left(\barint_{B_s(\mathbf{b})}
	|h-(h)_{B_s(\mathbf{b})}|^2\dx{x}\right)^{1/2}.
\]
We next claim that, for every \(\varepsilon>0\), there is
\(C_\varepsilon<\infty\) such that
\begin{equation}\label{eq:holomorphic-observation}
	\|h\|_{L^2(B_{1/2})}
	\le\varepsilon\|h\|_{L^2(B_1)}
	+C_\varepsilon\|h\|_{L^2(B_1\cap\{x_2<0\})}
\end{equation}
for every holomorphic \(h\in L^2(B_1)\). Otherwise, normalization
gives a sequence of holomorphic functions \(h_j\) satisfying
\[
	\|h_j\|_{L^2(B_{1/2})}=1,\qquad
	\varepsilon\|h_j\|_{L^2(B_1)}
	+j\|h_j\|_{L^2(B_1\cap\{x_2<0\})}<1.
\]
In particular, we get
\[
	\|h_j\|_{L^2(B_1)}<\varepsilon^{-1},\qquad
	\|h_j\|_{L^2(B_1\cap\{x_2<0\})}<j^{-1}.
\]

For \(0<t<1\) and \(z\in B_t\), applying the mean-value inequality for
subharmonic functions to \(|h_j|^2\) gives
\[
	|h_j(z)|^2
	\le\frac{4}{\pi(1-t)^2}
	\int_{B_{(1-t)/2}(z)}|h_j|^2\dx{x}
	\le\frac{4}{\pi(1-t)^2}\varepsilon^{-2}.
\]
Together with \eqref{eq:holomorphic-interior}, this gives local
uniform bounds for \(h_j\) and \(\nabla h_j\). By the Arzel\`a--Ascoli theorem,
a subsequence converges locally uniformly to a holomorphic function
\(h\). The bound on the lower half disk gives \(h=0\) there; hence it
follows from the standard identity theorem for holomorphic
functions that \(h\equiv0\) on \(B_1\). This contradicts
\[
	\|h\|_{L^2(B_{1/2})}
	=\lim_{j\to\infty}\|h_j\|_{L^2(B_{1/2})}=1.
\]
For \(\mathbf{a}=(a_1,0)\), apply \eqref{eq:holomorphic-observation} to
\(h(\mathbf{a}+r\zeta)\). Since
\[
	\|h(\mathbf{a}+r\zeta)\|_{L^2(B_t)}
	=r^{-1}\|h\|_{L^2(B_{tr}(\mathbf{a}))}
\]
for \(t=1\) and \(t=1/2\), and since \(a_2=0\), we obtain
\[
	\|h\|_{L^2(B_{r/2}(\mathbf{a}))}
	\le\varepsilon\|h\|_{L^2(B_r(\mathbf{a}))}
	+C_\varepsilon
	\|h\|_{L^2(B_r(\mathbf{a})\cap\{x_2<0\})},
\]
with constants independent of \(\mathbf{a}\) and \(r\).

Fix \(\mathbf{a}=(a_1,0)\), \(|a_1|\le1/4\), and \(0<r\le1/4\).
Since \(B_r(\mathbf{a})\subset B_{1/2}\), Lemma~\ref{lem:approximation}
gives \(\Phi,h_{B_r(\mathbf{a})}\in L^2(B_r(\mathbf{a}))\). Moreover,
\eqref{eq:approximation} gives
\begin{align*}
	\|\Phi-h_{B_r(\mathbf{a})}\|_{L^2(B_r(\mathbf{a}))}
	&\le C\bigl(r+\sup_{B_r(\mathbf{a})}|x_2|\bigr)
	|B_r(\mathbf{a})\cap\{u=0\}|^{1/2}\\
	&\le C(2r)(\pi r^2)^{1/2}\le Cr^2.
\end{align*}
Thus
\begin{equation}\label{eq:line-disk-error}
	\|\Phi-h_{B_r(\mathbf{a})}\|_{L^2(B_r(\mathbf{a}))}\le Cr^2.
\end{equation}
Since \(\Phi=0\) almost everywhere below the line, we deduce
\begin{equation}\label{eq:lower-half-observation}
	\begin{split}
		\|h_{B_r(\mathbf{a})}\|_{L^2(B_r(\mathbf{a})\cap\{x_2<0\})}
		&=\|h_{B_r(\mathbf{a})}-\Phi\|_{L^2(B_r(\mathbf{a})\cap\{x_2<0\})}\\
		&\le\|h_{B_r(\mathbf{a})}-\Phi\|_{L^2(B_r(\mathbf{a}))}\le Cr^2.
	\end{split}
\end{equation}

The rescaled form of \eqref{eq:holomorphic-observation}, together
with \eqref{eq:line-disk-error}--\eqref{eq:lower-half-observation},
now yields
\begin{align*}
	\|\Phi\|_{L^2(B_{r/2}(\mathbf{a}))}
	&\le\|\Phi-h_{B_r(\mathbf{a})}\|_{L^2(B_{r/2}(\mathbf{a}))}
	+\|h_{B_r(\mathbf{a})}\|_{L^2(B_{r/2}(\mathbf{a}))}\\
	&\le Cr^2+\varepsilon\|h_{B_r(\mathbf{a})}\|_{L^2(B_r(\mathbf{a}))}
	+C_\varepsilon
	\|h_{B_r(\mathbf{a})}\|_{L^2(B_r(\mathbf{a})\cap\{x_2<0\})}\\
	&\le Cr^2+\varepsilon\|\Phi\|_{L^2(B_r(\mathbf{a}))}
	+\varepsilon\|h_{B_r(\mathbf{a})}-\Phi\|_{L^2(B_r(\mathbf{a}))}
	+CC_\varepsilon r^2\\
	&\le\varepsilon\|\Phi\|_{L^2(B_r(\mathbf{a}))}
	+C(1+\varepsilon+C_\varepsilon)r^2.
\end{align*}
Thus, after adjusting \(C_\varepsilon\), we obtain
\begin{equation}\label{eq:hopf-decay-recursion}
	\|\Phi\|_{L^2(B_{r/2}(\mathbf{a}))}
	\le\varepsilon\|\Phi\|_{L^2(B_r(\mathbf{a}))}
	+C_\varepsilon r^2.
\end{equation}

Taking \(\varepsilon=1/16\) in
\eqref{eq:hopf-decay-recursion} and multiplying by \((r/2)^{-2}\),
we obtain
\[
	(r/2)^{-2}\|\Phi\|_{L^2(B_{r/2}(\mathbf{a}))}
	\le\frac14r^{-2}\|\Phi\|_{L^2(B_r(\mathbf{a}))}+C.
\]
At radius \(1/4\), \eqref{eq:initial-hopf-bound} gives
\[
	(1/4)^{-2}\|\Phi\|_{L^2(B_{1/4}(\mathbf{a}))}
	\le16\|\Phi\|_{L^2(B_{3/4})}
	\le C\left(1+\int_{B_1}|\nabla u|^2\,\dx{x}\right).
\]
Iteration at radii \(r=2^{-k-2}\) gives
\begin{align*}
	(2^{-j-2})^{-2}\|\Phi\|_{L^2(B_{2^{-j-2}}(\mathbf{a}))}
	&\le4^{-j}(1/4)^{-2}\|\Phi\|_{L^2(B_{1/4}(\mathbf{a}))}
	+C\sum_{k=0}^{j-1}4^{-k}\\
	&\le C\left(1+\int_{B_1}|\nabla u|^2\,\dx{x}\right)
	+\frac{C}{1-1/4}\\
	&\le C\left(1+\int_{B_1}|\nabla u|^2\,\dx{x}\right).
\end{align*}
For arbitrary \(r\in(0,1/4]\), choose \(j\ge0\) with
\(2^{-j-3}<r\le2^{-j-2}\). Then
\begin{align*}
	\|\Phi\|_{L^2(B_r(\mathbf{a}))}
	&\le\|\Phi\|_{L^2(B_{2^{-j-2}}(\mathbf{a}))}\\
	&\le C\left(1+\int_{B_1}|\nabla u|^2\,\dx{x}\right)
	(2^{-j-2})^2\\
	&\le4C\left(1+\int_{B_1}|\nabla u|^2\,\dx{x}\right)r^2.
\end{align*}
This proves \eqref{eq:decay}, uniformly in \(\mathbf{a}\) and \(r\).
\end{proof}

\begin{remark}
The interior estimate \eqref{eq:holomorphic-interior}, established
at the beginning of the proof of Lemma~\ref{lem:decay},
verifies the second approximation hypothesis
\eqref{eq:DHN-second} in Theorem~\ref{thm:DHN} when that theorem
is applied to the Hopf differential. In the abstract formulation
of De~Philippis--Hirsch--Nahon, \eqref{eq:DHN-second} is an
assumption on the local approximants; the approximants are not
required to be holomorphic or harmonic. In the present setting,
however, Lemma~\ref{lem:holomorphic-correction} constructs
approximants \(h_B\) that are holomorphic. Consequently,
\eqref{eq:holomorphic-interior}, and hence
\eqref{eq:DHN-second} with a universal constant \(K\), follows
naturally from standard interior estimates for holomorphic
functions. The remaining observation on the lower half disk, and the
decay iteration in Lemma~\ref{lem:decay}, serve a different purpose:
they establish the decay estimate \eqref{eq:decay}, from which the
height dependence is obtained.
\end{remark}

\begin{remark}
For the unweighted Bernoulli problem considered in \cite{DHN26}, 
the local \(L^\infty\)
bound for the Hopf differential
\cite[Lemma~22]{DHN26} (see also
\cite[Corollary~7]{DHN26} for the resulting Lipschitz estimate)
gives the general local bound
\[
	\|\Phi\|_{L^2(B_r)}\lesssim r,
\]
simply because \(|B_r|^{1/2}\simeq r\). In contrast,
\eqref{eq:decay} provides the improved order
\[
	\|\Phi\|_{L^2(B_r(\mathbf{a}))}\lesssim r^2
\]
for disks centered on the line \(\{x_2=0\}\). This additional
factor of \(r\) is crucial when applying
Theorem~\ref{thm:DHN}. Indeed, the approximation error gives
\(\Lambda\lesssim r\), while \eqref{eq:decay} yields
\[
	r^{-1}\|\Phi\|_{L^2(B_r(\mathbf{a}))}\lesssim r.
\]
Thus both terms in \eqref{eq:DHN-bound} are of order \(r\),
leading to \(\|\Phi\|_{L^\infty(B_{r/2}(\mathbf{a}))}\lesssim r\).
With only the \(O(r)\) bound for
\(\|\Phi\|_{L^2(B_r(\mathbf{a}))}\), the second term would merely be of
order \(1\), giving a gradient bound independent of the height rather
than the sharp estimate \(|\nabla u(x)|^2\lesssim x_2\).
\end{remark}

To obtain a supremum bound, we use the following scaled form of
\cite[Theorem~1]{DHN26} for maps with values in \(\C\simeq\R^2\).
For an integrable function on a disk \(B\), let
\((h)_B=|B|^{-1}\int_Bh\,\dx{x}\). The \emph{associated boundary}
of a measurable function \(W\) on \(B_r(\mathbf{a})\) is
\[
	\mathcal F_W=
	\left\{x\in B_r(\mathbf{a}):
	\begin{aligned}
		&\text{for every }0<\rho<r-|x-\mathbf{a}|,\\
		&|B_\rho(x)\cap\{W=0\}|\,|B_\rho(x)\cap\{W\ne0\}|>0
	\end{aligned}
	\right\}.
\]

\begin{theorem}[De Philippis--Hirsch--Nahon, Theorem 1 in \cite{DHN26}]\label{thm:DHN}
Let \(K\ge0\), \(\Lambda>0\), and
\(W\in L^2(B_r(\mathbf{a});\C)\). Assume that each disk
\(B_s(\mathbf{b})\subset B_r(\mathbf{a})\) admits
\(h_{B_s(\mathbf{b})}\in L^2(B_s(\mathbf{b});\C)\) satisfying
\begin{align}
	\int_{B_s(\mathbf{b})}|W-h_{B_s(\mathbf{b})}|^2\dx{x}
	&\le\Lambda^2|B_s(\mathbf{b})\cap\{W=0\}|,\label{eq:DHN-first}\\
	s\|\nabla h_{B_s(\mathbf{b})}\|_{L^\infty(B_{s/2}(\mathbf{b}))}
	&\le K\left(\barint_{B_s(\mathbf{b})}
	|h_{B_s(\mathbf{b})}-(h_{B_s(\mathbf{b})})_{B_s(\mathbf{b})}|^2\dx{x}\right)^{1/2}.
	\label{eq:DHN-second}
\end{align}
If \(\mathcal F_W\ne\varnothing\), then
\[
	|W(x)|\le C\left[
	\Lambda+
	\left(\frac{\operatorname{dist}(x,\mathcal F_W)}r\right)^\alpha
	r^{-1}\|W\|_{L^2(B_r(\mathbf{a}))}\right]
\]
for almost every \(x\in B_{r/2}(\mathbf{a})\), where \(C,\alpha>0\)
depend only on \(K\). Consequently,
\begin{equation}\label{eq:DHN-bound}
	\|W\|_{L^\infty(B_{r/2}(\mathbf{a}))}
	\le C\left(\Lambda+r^{-1}\|W\|_{L^2(B_r(\mathbf{a}))}\right).
\end{equation}
For holomorphic approximants, \eqref{eq:DHN-bound} remains valid
when \(\mathcal F_W=\varnothing\).
\end{theorem}

Since \(\operatorname{dist}(x,\mathcal F_W)/r\le3/2\) on
\(B_{r/2}(\mathbf{a})\), the value of \(\alpha\) is immaterial here.
Appendix~\ref{app:approximation-reduction} derives this scaled
form from \cite[Theorem~1]{DHN26} and proves the additional
assertion for holomorphic approximants when
\(\mathcal F_W=\varnothing\).

Applying Theorem~\ref{thm:DHN} to the Hopf differential \eqref{eq:hopf-definition} gives
the following estimate.

\begin{proposition}\label{prop:hopf-supremum}
Let \(u\in W^{1,2}(B_1)\) be nonnegative and vanish almost everywhere
on \(B_1\cap\{x_2\le0\}\), and suppose that \(u\) satisfies \eqref{eq:stationarity}
on \(\Omega=B_1\) for every \(\phi\in C_c^1(B_1;\R^2)\).
Set \(\Phi=(u_{x_1}-iu_{x_2})^2\). Then, for
\(\mathbf{a}=(a_1,0)\), \(|a_1|\le1/4\), and \(0<r\le1/4\),
\begin{equation}\label{eq:hopf-supremum}
	\|\Phi\|_{L^\infty(B_{r/2}(\mathbf{a}))}
	\le C\left(1+\int_{B_1}|\nabla u|^2\dx{x}\right)r,
\end{equation}
where \(C>0\) is universal.
\end{proposition}

\begin{proof}
Fix \(\mathbf{a}=(a_1,0)\), \(|a_1|\le1/4\), and \(0<r\le1/4\).
By Lemma~\ref{lem:approximation}, \(\Phi\in L^2(B_r(\mathbf{a}))\) and
each \(B_s(\mathbf{b})\subset B_r(\mathbf{a})\) admits a holomorphic approximant
\(h_{B_s(\mathbf{b})}\in L^2(B_s(\mathbf{b}))\). Since \(s\le r\) and
\(|x_2|\le r\) on this disk, \eqref{eq:approximation} yields
\begin{align*}
	\int_{B_s(\mathbf{b})}|\Phi-h_{B_s(\mathbf{b})}|^2\,\dx{x}
	&\le C\bigl(s+\sup_{B_s(\mathbf{b})}|x_2|\bigr)^2
	|B_s(\mathbf{b})\cap\{u=0\}|\\
	&\le Cr^2|B_s(\mathbf{b})\cap\{u=0\}|\\
	&\le Cr^2|B_s(\mathbf{b})\cap\{\Phi=0\}|.
\end{align*}
Here \(\nabla u=0\) almost everywhere on \(\{u=0\}\).
Thus \eqref{eq:DHN-first} holds with \(\Lambda=Cr\), while
\eqref{eq:DHN-second} follows from \eqref{eq:holomorphic-interior}
with universal \(K\). By \eqref{eq:DHN-bound} and \eqref{eq:decay}, we obtain
\begin{align*}
	\|\Phi\|_{L^\infty(B_{r/2}(\mathbf{a}))}
	&\le C\left(r+r^{-1}\|\Phi\|_{L^2(B_r(\mathbf{a}))}\right)\\
	&\le Cr+Cr^{-1}
	\left(1+\int_{B_1}|\nabla u|^2\,\dx{x}\right)r^2\\
	&\le Cr\left(1+\int_{B_1}|\nabla u|^2\,\dx{x}\right).
\end{align*}
\end{proof}

Proposition~\ref{prop:hopf-supremum} now gives the gradient bound
by choosing disks at the height of the point under consideration.

\begin{proof}[Proof of Theorem~\ref{thm:growth}]
We first work on \(B_1\). Let \(x=(x_1,x_2)\in B_{1/16}\) with \(x_2>0\) be a Lebesgue point
of \(\Phi\). Set \(\mathbf{a}=(x_1,0)\) and \(r=4x_2\). Then
\(|a_1|<1/4\), \(r<1/4\), and \(|x-\mathbf{a}|=r/4<r/2\).
The Lebesgue differentiation theorem gives
\[
	|\Phi(x)|
	=\lim_{\rho\to0^+}\barint_{B_\rho(x)}|\Phi|\,\dx{y}
	\le\|\Phi\|_{L^\infty(B_{r/2}(\mathbf{a}))}.
\]
Using \(r=4x_2\) and \(|\Phi|=|\nabla u|^2\) in \eqref{eq:hopf-supremum}, we obtain
\begin{equation}\label{eq:unit-gradient-conclusion}
	|\nabla u(x)|^2
	\le C\left(1+\int_{B_1}|\nabla u|^2\,\dx{y}\right)x_2
	\quad\hbox{for a.e. }x\in B_{1/16}\cap\{x_2>0\}.
\end{equation}
Since \(u=0\) almost everywhere below the line, we have
\(\nabla u\in L^\infty(B_{1/16})\). Hence \(u\) has a nonnegative
Lipschitz representative that vanishes on \(\{x_2\le0\}\).
Absolute continuity on almost every vertical segment and
\eqref{eq:unit-gradient-conclusion} give
\begin{align*}
	0\le u(x_1,x_2)
	&=\int_0^{x_2}u_{x_2}(x_1,t)\,\dx{t}\\
	&\le\int_0^{x_2}|\nabla u(x_1,t)|\,\dx{t}\\
	&\le C\left(1+\int_{B_1}|\nabla u|^2\,\dx{y}\right)^{1/2}
	\int_0^{x_2}t^{1/2}\,\dx{t}\\
	&=\frac{2C}{3}
	\left(1+\int_{B_1}|\nabla u|^2\,\dx{y}\right)^{1/2}
	x_2^{3/2}.
\end{align*}
By continuity, the estimate holds throughout the disk. For
\(\mathbf{b}\in B_{1/16}\cap\{x_2=0\}\), it also gives
\[
	\frac{|u(x)-u(\mathbf{b})|}{|x-\mathbf{b}|}
	\le C\left(1+\int_{B_1}|\nabla u|^2\,\dx{y}\right)^{1/2}
	|x-\mathbf{b}|^{1/2}\longrightarrow0.
\]
Thus \(Du(\mathbf{b})=0\). By Definition~\ref{def:variational}, \(u\in C^2(\Omega\cap\{u>0\})\),
so \eqref{eq:unit-gradient-conclusion} holds pointwise in the positive phase.

For the disk \(B_R(x^0)\) with \(x_2^0=0\), set
\(u_R(x)=R^{-3/2}u(x^0+Rx)\). We test \eqref{eq:stationarity}
with \(\phi_R(y)=R\phi((y-x^0)/R)\), where
\(\phi\in C_c^1(B_1;\R^2)\). Writing \(y=x^0+Rx\), we have
\[
	\begin{aligned}
		\nabla u_R(x)&=R^{-1/2}\nabla u(y),&
		\nabla\phi_R(y)&=\nabla\phi(x),\\
		(\phi_R)_2(y)&=R\phi_2(x),& y_2&=Rx_2.
	\end{aligned}
\]
Since \(\dx{y}=R^2\dx{x}\) and
\(\chi_{\{u>0\}}(y)=\chi_{\{u_R>0\}}(x)\), the identity becomes
\begin{align*}
	0={}&\int_{B_R(x^0)}\bigl[
	(|\nabla u|^2+y_2\chi_{\{u>0\}})\diver\phi_R
	-2\nabla uD\phi_R\nabla u
	+\chi_{\{u>0\}}(\phi_R)_2\bigr]\dx{y}\\
	={}&R^3\int_{B_1}\bigl[
	(|\nabla u_R|^2+x_2\chi_{\{u_R>0\}})\diver\phi
	-2\nabla u_RD\phi\nabla u_R
	+\chi_{\{u_R>0\}}\phi_2\bigr]\dx{x}.
\end{align*}
Dividing by \(R^3\) proves \eqref{eq:stationarity} for \(u_R\).
Moreover,
\[
	\int_{B_1}|\nabla u_R|^2\,\dx{x}
	=R^{-1}\int_{B_1}|\nabla u(x^0+Rx)|^2\,\dx{x}
	=R^{-3}\int_{B_R(x^0)}|\nabla u|^2\,\dx{y}.
\]
Applying the unit disk estimates to \(u_R\) and scaling back gives
\eqref{eq:main-gradient}--\eqref{eq:main-function}, together with
local Lipschitz continuity and differentiability on the line \(\{x_2=0\}\). A finite covering makes the growth constant uniform near
compact portions of the line \(\{x_2=0\}\) in \(\Omega\).

\end{proof}

\begin{remark}
For completeness we record that the arguments above also give the
local bound \eqref{eq:local-growth-intro} at points away from the
line \(\{x_2=0\}\), with a constant that may depend on the local
disk and the local energy. Indeed, let
\(B_s(\mathbf b)\Subset\Omega\cap\{x_2>0\}\) and set
\(H=\sup_{B_s(\mathbf b)}|x_2|+2s\). Since the derivation of
\eqref{eq:hopf} is local, it holds on every disk
\(B_t(\mathbf c)\subset B_s(\mathbf b)\); the correction
\eqref{eq:holomorphic-correction} and \eqref{eq:correction-bound}
then verify the first hypothesis of
Theorem~\ref{thm:DHN} with \(\Lambda=H\), while
\eqref{eq:holomorphic-interior} gives the second one with a
universal \(K\). Applying the mean-value argument of
Lemma~\ref{lem:integrability-height} on
\(B_s(\mathbf b)\) yields
\[
\|\Phi\|_{L^2(B_{s/2}(\mathbf b))}
\le C\Bigl(s^{-1}\int_{B_s(\mathbf b)}|\nabla u|^2\,\dx{x}+sH\Bigr)
\]
with \(C\) universal. Applying \eqref{eq:DHN-bound} on
\(B_{s/2}(\mathbf b)\), we obtain
\[
\|\Phi\|_{L^\infty(B_{s/4}(\mathbf b))}
\le C\Bigl(H+s^{-2}\int_{B_s(\mathbf b)}|\nabla u|^2\,\dx{x}\Bigr).
\]
Since \(B_{s/4}(\mathbf b)\Subset\{x_2>0\}\), division by
\(\inf_{B_{s/4}(\mathbf b)}x_2>0\) gives
\(|\nabla u|^2\le Cx_2\) almost everywhere on
\(B_{s/4}(\mathbf b)\). Below the line \(\{x_2=0\}\) we have
\(u=0\). Hence \eqref{eq:local-growth-intro} holds locally in
\(\Omega\).
\end{remark}

The following lemma verifies that the Stokes profile belongs to the
variational class and proves the optimality of the exponents.

\begin{lemma}
Let \(u=\US\) be the Stokes profile in \eqref{eq:stokes-profile}.
Then \(u\) is a variational solution of \eqref{eq.problem} on every
bounded domain in \(\R^2\), and
\[
	|\nabla u(0,x_2)|^2=\frac{x_2}{2},\qquad
	u(0,x_2)=\frac{\sqrt2}{3}x_2^{3/2}\quad(x_2>0).
\]
Consequently, the powers of \(x_2\) in
\eqref{eq:main-gradient} and \eqref{eq:main-function} are optimal.
\end{lemma}

\begin{proof}
Differentiating \eqref{eq:stokes-profile} in the positive sector
\(\{u>0\}=\{\pi/6<\theta<5\pi/6\}\) gives
\[
	\partial_\rho u
	=\frac{\sqrt2}{2}\rho^{1/2}
	\cos\!\left(\frac32\left(\theta-\frac\pi2\right)\right),
	\qquad
	\rho^{-1}\partial_\theta u
	=-\frac{\sqrt2}{2}\rho^{1/2}
	\sin\!\left(\frac32\left(\theta-\frac\pi2\right)\right).
\]
Thus \(u\) is harmonic there and \(|\nabla u|^2=\rho/2\).
The zero extension is continuous and belongs to
\(W^{1,2}_{\loc}(\R^2)\). On each boundary ray,
\(x_2=\rho/2=|\nabla u|^2\) holds.

Fix \(\phi\in C_c^1(\R^2;\R^2)\) and choose \(R>0\) with
\(\operatorname{supp}\phi\Subset B_R\). Integration by parts
in \((B_R\setminus\overline{B_\varepsilon})\cap\{u>0\}\)
gives
\[
	\begin{split}
		&\int_{B_R\setminus B_\varepsilon}
		\left[(|\nabla u|^2+x_2\chi_{\{u>0\}})\diver\phi
		-2\nabla uD\phi\nabla u+\chi_{\{u>0\}}\phi_2\right]\dx{x}\\
		&\quad=\int_{\partial B_\varepsilon\cap\{u>0\}}
		\left[(|\nabla u|^2+x_2)\phi\cdot\nu
		-2(\nabla u\cdot\nu)(\nabla u\cdot\phi)\right]\dH{1},
	\end{split}
\]
where \(\nu\) points into \(B_\varepsilon\). The term on the outer
circle vanishes because \(\phi=0\) there, and the terms on the
rays vanish because their integrands are
\((x_2-|\nabla u|^2)\phi\cdot\nu\).
On the inner circle we have \(|\nabla u|^2+x_2\le C\varepsilon\),
so the flux through the inner arc is \(O(\varepsilon^2)\).
The omitted quadratic and weighted terms over \(B_\varepsilon\)
are \(O(\varepsilon^3)\), while
\(\int_{B_\varepsilon}\chi_{\{u>0\}}\phi_2\dx{x}
=O(\varepsilon^2)\). Letting \(\varepsilon\to0^+\) proves
\eqref{eq:stationarity}. The identities on the vertical axis
follow by setting \(\theta=\pi/2\), and show that neither
exponent can be increased.
\end{proof}

\FloatBarrier

\section{The Stokes conjecture for variational solutions}

Throughout this section, \(u\) is a variational solution of \eqref{eq.problem}. Theorem~\ref{thm:growth}
provides the local gradient bound. We first establish the measure
and perimeter estimates of \cite{DPY26} in the two-dimensional irrotational
case, and then prove the three assertions of
Theorem~\ref{thm:stokes} in order.
Appendix~\ref{app:hypothesis-audit} summarizes the role of the
variational, topological, and free boundary assumptions in
applying the classification results of \cite{VW11} and proving
the final regularity assertion.

\begin{lemma}[Boundary measure and weighted perimeter]
\label{lem:variational-measure-perimeter}
Let \(u\) be a variational solution of \eqref{eq.problem} in
\(B_R(x^0)\), where \(x^0=(x_1^0,0)\), and assume that
\[
	|\nabla u|^2\le C_0x_2^+
	\quad\text{a.e. in }B_R(x^0)
\]
for some \(C_0>0\). Then the following hold.
\begin{itemize}
\item \(u\) is harmonic in \(\{u>0\}\), and \(\mu=\Delta u\) is a nonnegative Radon
measure supported on the free boundary.
\item The positive phase satisfies
\begin{equation}\label{eq:variational-weighted-perimeter}
	\chi_{\{u>0\}}\in BV_{\loc}(B_R(x^0)\cap\{x_2>0\}),
	\qquad \sqrt{x_2}\,|D\chi_{\{u>0\}}|\le\sqrt{C_0}\,\mu
\end{equation}
as measures in \(B_R(x^0)\cap\{x_2>0\}\), and
\begin{equation}\label{eq:variational-boundary-mass}
	\mu(B_r(x^0))\le C\sqrt{C_0}\,r^{3/2},
	\qquad 0<2r<R.
\end{equation}
\end{itemize}
\end{lemma}

\begin{proof}
We follow \cite[Lemmas~3.4--3.5]{DPY26}. After translation,
assume \(x^0=0\). The gradient bound gives local Lipschitz continuity,
and \(\nabla u=0\) almost everywhere on \(\{u=0\}\).
For test functions supported in \(\{u>0\}\), the gravitational
contribution to \eqref{eq:stationarity} vanishes because
\[
	\int\bigl(x_2\diver\phi+\phi_2\bigr)\dx{x}
	=\int\diver(x_2\phi)\dx{x}=0.
\]
Integration by parts in the remaining terms gives
\((\Delta u)\nabla u=0\). If \(\Delta u\ne0\) at a
point, continuity would force \(\nabla u=0\) in a neighborhood
of that point, a contradiction. Hence \(\Delta u=0\) in \(\{u>0\}\).

Choose a smooth nondecreasing \(\eta:\R\to[0,1]\) with
\(\eta=0\) on \(({-}\infty,1]\) and \(\eta=1\) on
\([2,\infty)\). Testing harmonicity with \(\eta(u/\varepsilon)\zeta\),
\(\zeta\in C_c^1(B_R)\), gives
\begin{equation}\label{eq:cusp-scalar-cutoff}
	\int_{B_R}\zeta\,
	\frac{\eta'(u/\varepsilon)}{\varepsilon}
	|\nabla u|^2\dx{x}
	=-\int_{B_R}\eta(u/\varepsilon)\nabla u\cdot\nabla\zeta\dx{x}.
\end{equation}
The measures
\[
	\mu_\varepsilon=
	\frac{\eta'(u/\varepsilon)}{\varepsilon}|\nabla u|^2\dx{x}
\]
are nonnegative. Dominated convergence on the right of
\eqref{eq:cusp-scalar-cutoff}, together with the vanishing of
\(\nabla u\) on the zero set, gives
\[
	\lim_{\varepsilon\to0^+}\int\zeta\,\dx{\mu_\varepsilon}
	=-\int\nabla u\cdot\nabla\zeta\dx{x}
	=\langle\Delta u,\zeta\rangle.
\]
Thus \(\mu=\Delta u\) is a nonnegative Radon measure, supported
on the free boundary by harmonicity. Taking a nonnegative cutoff function
equal to one on a given compact set in \eqref{eq:cusp-scalar-cutoff}
bounds \(\mu_\varepsilon\) there. Consequently,
\(\mu_\varepsilon\rightharpoonup\mu\) locally as measures.

To estimate \(\mu(B_r)\), take \(0\le\zeta\le1\) with
\(\zeta=1\) on \(B_r\), support in \(B_{2r}\), and
\(|\nabla\zeta|\le C/r\). Then
\[
	\mu(B_r)\le\int\zeta\,\dx{\mu}
	=-\int\nabla u\cdot\nabla\zeta\dx{x}
	\le\frac{C\sqrt{C_0}}r\int_{B_{2r}}\sqrt{x_2^+}\dx{x}
	\le C\sqrt{C_0}\,r^{3/2},
\]
which proves \eqref{eq:variational-boundary-mass}.

We next recover the perimeter from stationarity. Fix
\(\phi\in C_c^1(B_R\cap\{x_2>0\};\R^2)\).
Since \(\Delta u=0\) in the positive phase, integration by parts
against \(\eta(u/\varepsilon)\phi\) gives
\begin{equation}\label{eq:cusp-stationarity-cutoff}
	\int\eta(u/\varepsilon)
	\bigl(|\nabla u|^2\operatorname{div}\phi
	-2\nabla uD\phi\nabla u\bigr)\dx{x}
	=\int\nabla u\cdot\phi\,\dx{\mu_\varepsilon}.
\end{equation}
The gradient bound gives
\[
	\left|\int\nabla u\cdot\phi\,\dx{\mu_\varepsilon}\right|
	\le\sqrt{C_0}\int\sqrt{x_2}|\phi|\,\dx{\mu_\varepsilon}.
\]
Passing to the limit by weak convergence of the measures and the
dominated convergence theorem in \eqref{eq:cusp-stationarity-cutoff},
and using \eqref{eq:stationarity}, we obtain
\[
	\left|\int\chi_{\{u>0\}}
	\operatorname{div}(x_2\phi)\dx{x}\right|
	\le\sqrt{C_0}\int\sqrt{x_2}|\phi|\,\dx{\mu}.
\]
For \(\zeta\in C_c^1(B_R\cap\{x_2>0\};\R^2)\), substitute
\(\phi=\zeta/x_2\) to obtain
\[
	\left|\int\chi_{\{u>0\}}\operatorname{div}\zeta\dx{x}\right|
	\le\sqrt{C_0}\int\frac{|\zeta|}{\sqrt{x_2}}\,\dx{\mu}.
\]
Since the right hand side is locally finite above the line, the
dual characterization of \(BV\) gives
\(|D\chi_{\{u>0\}}|\le\sqrt{C_0}x_2^{-1/2}\mu\), which is
\eqref{eq:variational-weighted-perimeter}.
\end{proof}

Fix \(x^0\in\mathcal S^u\) and choose \(R>0\) such that
\(B_{32R}(x^0)\Subset\Omega\). Theorem~\ref{thm:growth}, applied
on this outer ball, gives
\[
	\begin{gathered}
		|\nabla u|^2\le C_0x_2^+
		\quad\text{a.e. in }B_{2R}(x^0),\\
		C_0=C\left(1+(32R)^{-3}
		\int_{B_{32R}(x^0)}|\nabla u|^2\dx{x}\right).
	\end{gathered}
\]
We use this bound on the smaller working balls below. It supplies
the growth hypothesis in the compactness and density results of
\cite{VW11}. Write
\[
	u_r(x)=r^{-3/2}u(x^0+rx).
\]
Every sequence \(r\to0^+\) has a subsequence along which
\(u_r\) converges locally uniformly and strongly in
\(W^{1,2}_{\loc}(\R^2)\), and each limit is homogeneous of degree
\(3/2\); see \cite[Lemma~4.1]{VW11}. We record also the energy identity
\begin{equation}\label{eq:energy-conclusion}
	\int\zeta|\nabla u|^2\dx{x}
	=-\int u\nabla u\cdot\nabla\zeta\dx{x},\qquad \zeta\in C_c^1(\Omega),
\end{equation}
which follows by testing harmonicity with \((u-\varepsilon)^+\zeta\)
and passing to the limit as \(\varepsilon\to0^+\). It also yields
strong convergence when the centers of the rescalings vary along
the line \(\{x_2=0\}\). Applying Lemma~\ref{lem:variational-measure-perimeter}
on \(B_{2R}(x^0)\) gives
\(\Delta u\ge0\) and
\begin{equation}\label{eq:weighted-perimeter-conclusion}
	\begin{gathered}
		\chi_{\{u>0\}}\in BV_{\loc}(B_R(x^0)\cap\{x_2>0\}),\\
		\int_{B_r(x^0)\cap\{x_2>0\}}\sqrt{x_2}\,
		\dx{|D\chi_{\{u>0\}}|}\le CC_0r^{3/2},\qquad 0<r<R.
	\end{gathered}
\end{equation}
Here \(C\) is universal and \(C_0\) is the gradient constant fixed
above. This verifies the weighted perimeter hypothesis of
\cite[Proposition~4.7]{VW11}, which gives the density trichotomy
\eqref{eq:boundary-density-conclusion} below.

The boundary Weiss functional is
\begin{equation}\label{eq:boundary-Weiss-conclusion}
	\begin{split}
		\Phi^{\mathrm{bound}}_{x^0,u}(r)
		={}&r^{-3}\int_{B_r(x^0)}
		\bigl(|\nabla u|^2+x_2\chi_{\{u>0\}}\bigr)\dx{x}\\
		&-\frac32r^{-4}\int_{\partial B_r(x^0)}u^2\dH{1}.
	\end{split}
\end{equation}
For a center \(x^0\) on the line \(\{x_2=0\}\), this functional is
locally absolutely continuous in the admissible positive radii.
The monotonicity formula \cite[Theorem~3.5(ii)]{VW11} gives
\begin{equation}\label{eq:boundary-Weiss-derivative-conclusion}
	\frac{\dx{}}{\dx{r}}\Phi^{\mathrm{bound}}_{x^0,u}(r)
	=2r^{-3}\int_{\partial B_r(x^0)}
	\left(\nabla u\cdot\nu-\frac{3u}{2r}\right)^2\dH{1}
	\ge0,
\end{equation}
for almost every such \(r\). In particular, if \(0<s<t\) and
\(\overline{B_t(x^0)}\subset\Omega\), then
\begin{equation}\label{eq:boundary-Weiss-integrated-conclusion}
	\begin{split}
		\Phi^{\mathrm{bound}}_{x^0,u}(t)-\Phi^{\mathrm{bound}}_{x^0,u}(s)
		=2\int_s^t\rho^{-3}\int_{\partial B_\rho(x^0)}
		\left(\nabla u\cdot\nu-\frac{3u}{2\rho}\right)^2
		\dH{1}\,\dx{\rho}.
	\end{split}
\end{equation}
The integrated form will be useful when passing to limits.
The possible boundary densities are given by
\cite[Lemmas~4.1--4.2 and Proposition~4.7]{VW11}:
\begin{equation}\label{eq:boundary-density-conclusion}
	\Phi^{\mathrm{bound}}_{x^0,u}(0+)
	=\lim_{r\to0^+}r^{-3}
	\int_{B_r(x^0)}x_2\chi_{\{u>0\}}\dx{x}
	\in\left\{0,\frac{\sqrt3}{3},\frac23\right\}.
\end{equation}

We use the variational cusp theorem of \cite[Theorem~I]{DPY26}.
For completeness, the following proof gives its two-dimensional
irrotational form and keeps the dependence on the nonsharp
gradient constant explicit.

\begin{proposition}[Exclusion of cusp density]
\label{prop:no-cusps-variational}
Let \(u\) be a variational solution of \eqref{eq.problem}, and let
\(x^0\in\mathcal S^u\). Assume that \(B_R(x^0)\Subset\Omega\) and
\[
	|\nabla u|^2\le C_0x_2^+
	\quad\text{a.e. in }B_R(x^0)
\]
for some \(C_0>0\). Then, for every \(0<r<R\),
\begin{equation}\label{eq:cusp-density-lower-bound}
	r^{-3}\int_{B_r(x^0)}x_2\chi_{\{u>0\}}\dx{x}
	\ge\frac{2}{\pi(1+C_0)^5}
	\left(\int_{B_1\cap\{x_2>0\}}\sqrt{x_2}\dx{x}\right)^2>0.
\end{equation}
Consequently, by Theorem~\ref{thm:growth}, the boundary Weiss
density at every stagnation point of a variational solution is positive.
\end{proposition}

\begin{proof}
Assume \(x^0=0\) and set
\[
	\begin{split}
		V(r)&=\int_{B_r}\sqrt{x_2^+}\chi_{\{u>0\}}\dx{x},\qquad
		I_2=\int_{B_1\cap\{x_2>0\}}\sqrt{x_2}\dx{x},\\
		P(r)&=\int_{B_r\cap\{x_2>0\}}\sqrt{x_2}
		\dx{|D\chi_{\{u>0\}}|}.
	\end{split}
\]
Lemma~\ref{lem:variational-measure-perimeter} gives
\(\mu=\Delta u\ge0\) and, by exhaustion above the line,
\begin{equation}\label{eq:cusp-perimeter-mass}
	P(r)\le\sqrt{C_0}\,\mu(B_r)<\infty\qquad(0<r<R).
\end{equation}
The coarea formula makes \(V\) absolutely continuous and gives
\begin{equation}\label{eq:cusp-volume-derivative}
	V'(r)=\int_{\partial B_r}\sqrt{x_2^+}
	\chi_{\{u>0\}}\dH{1}
\end{equation}
for almost every \(r\in(0,R)\).

Test \(\Delta u=\mu\) with a radial cutoff function equal to one on
\(B_r\), linear on \(B_{r+h}\setminus B_r\), and zero outside
\(B_{r+h}\). We then deduce from the smooth approximation, the coarea
formula, and the passage to the limit \(h\to0^+\) that
\[
	\mu(B_r)=\int_{\partial B_r}\nabla u\cdot\nu\dH{1}
	\le\sqrt{C_0}\int_{\partial B_r}\sqrt{x_2^+}
	\chi_{\{u>0\}}\dH{1}
	=\sqrt{C_0}V'(r)
\]
for almost every \(r\), for which also \(\mu(\partial B_r)=0\).
Here we used \(\nabla u=0\) on the zero set and
\eqref{eq:cusp-volume-derivative}. Together with
\eqref{eq:cusp-perimeter-mass}, this gives
\begin{equation}\label{eq:cusp-flux-perimeter}
	P(r)\le C_0V'(r)\qquad\text{for almost every }r\in(0,R).
\end{equation}

We next use the sharp isoperimetric inequality with weight
\(x_2^{1/2}\):
\begin{equation}\label{eq:cusp-weighted-isoperimetry}
	\int_{\{x_2>0\}}\sqrt{x_2}\,\dx{|D\chi_A|}
	\ge\frac52 I_2^{2/5}
	\left(\int_A\sqrt{x_2}\dx{x}\right)^{3/5};
\end{equation}
see \cite[Theorem~1.4]{CR13} and
\cite[Theorem~1.3]{CRS16} for the measurable set formulation.
Fix a radius for which the preceding identities hold and the
\(BV\) trace on the sphere agrees with its slice. By slicing on
a countable exhaustion of \(B_R\cap\{x_2>0\}\), these radii
have full measure.

For \(0<\varepsilon<r\), put
\(A_\varepsilon=B_r\cap\{u>0\}\cap\{x_2>\varepsilon\}\).
Since \(\chi_{\{u>0\}}\) has finite variation near the closed
cap, restriction and extension by zero give finite perimeter of
\(A_\varepsilon\) in \(\R^2\). Choose \(\varepsilon\to0^+\)
along heights at which the horizontal trace agrees with the slice.
The \(BV\) restriction formula
\cite[Section 5.3, Theorem 5.6 and Section 5.4, Theorem 5.8]{EG15}
(the truncated cap \(B_r\cap\{x_2>\varepsilon\}\) has Lipschitz
boundary, so both results apply) gives
\[
	\begin{split}
		\int_{\{x_2>0\}}\sqrt{x_2}\,\dx{|D\chi_{A_\varepsilon}|}
		\le{}&\int_{B_r\cap\{x_2>\varepsilon\}}\sqrt{x_2}
		\dx{|D\chi_{\{u>0\}}|}\\
		&+\int_{\partial B_r\cap\{x_2>\varepsilon\}}\sqrt{x_2}
		\chi_{\{u>0\}}\dH{1}\\
		&+\sqrt\varepsilon\int_{B_r\cap\{x_2=\varepsilon\}}
		\chi_{\{u>0\}}\dH{1}\\
		\le{}&P(r)+V'(r)+2r\sqrt\varepsilon.
	\end{split}
\]
The horizontal contribution is at most \(2r\sqrt\varepsilon\),
and the weighted volumes increase to \(V(r)\). Applying
\eqref{eq:cusp-weighted-isoperimetry} to \(A_\varepsilon\) and
letting \(\varepsilon\to0^+\) therefore gives
\begin{equation}\label{eq:cusp-radial-isoperimetry}
	P(r)+V'(r)\ge\frac52 I_2^{2/5}V^{3/5}(r)
\end{equation}
for almost every \(r<R\).

Combining \eqref{eq:cusp-flux-perimeter} and
\eqref{eq:cusp-radial-isoperimetry}, we obtain
\[
	V'(r)\ge\frac{5}{2(1+C_0)}I_2^{2/5}V^{3/5}(r).
\]
Since \(0\in\partial\{u>0\}\) and \(u\) is continuous,
\(V(r)>0\) for every \(r>0\). Thus, on any interval
\([\rho,r]\subset(0,R)\), the chain rule gives
\[
	\frac{\dx{}}{\dx{s}}V^{2/5}(s)
	=\frac25 V^{-3/5}(s)V'(s)
	\ge\frac{I_2^{2/5}}{1+C_0}
\]
for almost every \(s\). Integration and then \(\rho\to0^+\)
are justified by \(0\le V(\rho)\le I_2\rho^{5/2}\), and yield
\begin{equation}\label{eq:cusp-volume-lower-bound}
	V(r)\ge\frac{I_2}{(1+C_0)^{5/2}}\,r^{5/2}.
\end{equation}
For every \(r<R\), the Cauchy--Schwarz inequality gives
\[
	V^2(r)\le |B_r\cap\{u>0\}|
	\int_{B_r}x_2\chi_{\{u>0\}}\dx{x}
	\le\frac{\pi r^2}{2}
	\int_{B_r}x_2\chi_{\{u>0\}}\dx{x}.
\]
Combining this with \eqref{eq:cusp-volume-lower-bound} proves
\eqref{eq:cusp-density-lower-bound}. The last assertion follows
from \eqref{eq:boundary-density-conclusion} and
Theorem~\ref{thm:growth}.
\end{proof}

\begin{proof}[The two alternatives in Theorem~\ref{thm:stokes}\textup{(i)}]
Fix \(x^0\in\mathcal S^u\). By Theorem~\ref{thm:growth} and
Proposition~\ref{prop:no-cusps-variational}, the density in
\eqref{eq:boundary-density-conclusion} belongs to
\(\{1/\sqrt3,2/3\}\).

At density \(1/\sqrt3\), the classification of homogeneous limits
in \cite[Proposition~4.7]{VW11} identifies every blow-up limit
with \(\US\). At density \(2/3\), every blow-up limit is zero.
The compactness in \cite[Lemma~4.1]{VW11} then gives convergence
of the whole family \(u_r\), locally uniformly and strongly in
\(W^{1,2}_{\loc}(\R^2)\), since every subsequential limit is the
same.

For phase convergence, \eqref{eq:weighted-perimeter-conclusion}
gives local \(BV\) compactness in every region
\(\{x_2>\tau\}\), \(\tau>0\). From any sequence \(r_j\to0^+\),
we may therefore extract a diagonal subsequence such that, using
the vanishing area of strips along the line,
\[
	\chi_{\{u_{r_j}>0\}}\longrightarrow\chi_0
	\quad\text{in }L^1_{\loc}(\R^2),\qquad
	\chi_{\{u_0>0\}}\le\chi_0\le\chi_{\{x_2>0\}},
\]
where \(\chi_0\) is a characteristic function. The lower
inequality follows from local uniform convergence wherever
\(u_0>0\); the upper inequality follows from the fixed zero
phase below the line. For each \(R>0\), convergence of the
phase indicators and \eqref{eq:boundary-density-conclusion} yield
\[
	\int_{B_R}x_2\chi_0\dx{x}
	=R^3\Phi^{\mathrm{bound}}_{x^0,u}(0+).
\]
If \(u_0=\US\), the right hand side is
\(R^3/\sqrt3=\int_{B_R}x_2\chi_{\{\US>0\}}\dx{x}\).
The nonnegative difference
\(x_2(\chi_0-\chi_{\{\US>0\}})\) thus has zero integral
on every disk, so \(\chi_0=\chi_{\{\US>0\}}\) almost
everywhere. In the full-density case, i.e.\ when \eqref{eq:main-flat-blowup}
holds, the right hand side is
\(2R^3/3=\int_{B_R}x_2^+\dx{x}\), and the upper bound on
\(\chi_0\) similarly gives
\(\chi_0=\chi_{\{x_2>0\}}\) almost everywhere.
Uniqueness of the phase limit gives convergence of the whole
family and completes the proof.
\end{proof}

We now assume the finite component and interface hypotheses in
Theorem~\ref{thm:stokes}\textup{(ii)}. The singular set above
the line \(\{x_2=0\}\) has \(\mathcal H^1\) measure zero. The
one sided Bernoulli trace on the regular patches gives the
lower bound for the boundary measure.

Define the full-density set by
\[
	\Sigma^u=\left\{x^0\in \mathcal S^u:
	\lim_{r\to0^+}r^{-3}\int_{B_r(x^0)}
	x_2\chi_{\{u>0\}}\dx{x}=\frac23\right\}.
\]
For use in the finite component argument, we record the two
frequency conclusions at a point \(x^0\in\Sigma^u\).
There is an integer \(N(x^0)\ge2\) such that, for
\(0\le\beta<N(x^0)\),
\begin{equation}\label{eq:subfrequency-vanishing}
	\frac{u\bigl((x_1^0,0)+r(x_1,x_2)\bigr)}{r^\beta}\longrightarrow0
\end{equation}
strongly in \(L^2_{\loc}(\R^2)\) in the spatial variables
\((x_1,x_2)\) as \(r\to0^+\). With
\(x=\rho(\cos\theta,\sin\theta)\), \(0\le\theta<2\pi\), we obtain
\begin{equation}\label{eq:degenerate-profile}
	\frac{u(x^0+rx)}
	{\left(r^{-1}\int_{\partial B_r(x^0)}u^2\dH{1}\right)^{1/2}}
	\longrightarrow
	\begin{cases}
		\displaystyle\sqrt{\frac2\pi}\,\rho^{N(x^0)}
		|\sin(N(x^0)\theta)|,&0<\theta<\pi,\\[2pt]
		0,&\text{otherwise},
	\end{cases}
\end{equation}
strongly in \(W^{1,2}_{\loc}(B_1\setminus\{0\})\) and weakly in
\(W^{1,2}(B_1)\).

Part~\textup{(i)} gives the Stokes limit and convergence of
the phases outside \(\Sigma^u\). The density and frequency results in
\cite[Lemma~4.2(iv), Remark~6.2 and Theorems~10.1--10.2]{VW11}
give the relative closedness and local finiteness of \(\Sigma^u\), together with
the normalized profile \eqref{eq:degenerate-profile}. Its normalization uses
\(\int_0^\pi\sin^2(N\theta)\,\dx{\theta}=\pi/2\).\footnote{In
\cite[Theorem~10.1]{VW11}, the denominator is
\(\bigl(\int_0^\pi\sin^2(N\theta)\,\dx{\theta}\bigr)^{1/2}\),
which gives the coefficient \(\sqrt{2/\pi}\) used here.}
The assertion
\eqref{eq:subfrequency-vanishing} follows from
\cite[Lemma~7.5(iv)]{VW11}: by \cite[Theorem~10.1]{VW11}, the
limiting frequency is an integer \(N(x^0)\ge2\).
These facts supply the blow-up and frequency information needed
for the exclusion of full density below.

\begin{proposition}\label{prop:no-full-density}
Let \(u\) be a variational solution of \eqref{eq.problem}, and assume
the finite component and free boundary hypotheses in
Theorem~\ref{thm:stokes}\textup{(ii)}. Then
\(\Sigma^u=\varnothing\).
\end{proposition}

\begin{proof}
Suppose by contradiction that \(x^0\in\Sigma^u\), and translate
\(x^0\) to the origin. Fix an open neighborhood
\(U\Subset\Omega\) in which \(\{u=0\}\cap U\) has finitely
many connected components, and let \(C\) be the component
containing \(0\). The other components are relatively closed
in \(U\), so a sufficiently small ball about \(0\) meets none
of them.

Choose \(r_j\to0^+\) and set
\[
	u_j(x)=r_j^{-3/2}u(r_jx),\qquad
	v_j(x)=\frac{u(r_jx)}
	{\left(r_j^{-1}\int_{\partial B_{r_j}}u^2\dH{1}\right)^{1/2}}.
\]
Their positive and zero sets agree. By part~\textup{(i)},
\(u_j\) converges to 0 locally uniformly and
\(\chi_{\{u_j>0\}}\) converges to \(\chi_{\{x_2>0\}}\) in
\(L^1_{\loc}(\R^2)\). By \eqref{eq:degenerate-profile},
\(v_j\) converges in \(L^2_{\loc}(B_1\setminus\{0\})\) to
\[
	v(\rho,\theta)=\sqrt{2/\pi}\,\rho^N|\sin(N\theta)|,
	\qquad 0<\theta<\pi,\qquad N=N(0)\ge2.
\]
Choose a point \(z\) on an interior nodal ray of \(v\), and
\(\delta>0\) such that
\(\overline{B_{4\delta}(z)}\subset B_1\cap\{x_2>0\}\)
and this disk meets no other nodal ray. The function \(v\) is
not harmonic in \(B_\delta(z)\): there it is the absolute value
of a harmonic function whose normal derivative on the nodal ray
is nonzero.

For all sufficiently large \(j\),
\(\Gamma_j:=\partial\{u_j>0\}\) meets \(B_\delta(z)\).
Otherwise, along a subsequence, the connected disk contains no
free boundary point, so \(v_j\) is harmonic throughout it,
either in the positive phase or as the zero function. Passing
to the limit in distributions would make \(v\) harmonic there,
a contradiction. Choose \(\mathbf{p}_j\in\Gamma_j\cap B_\delta(z)\).
For large \(j\), the component reduction above implies that
\(\mathbf{p}_j\) and \(0\) belong to the connected set
\(C_j:=r_j^{-1}C\). The continuous map
\(x\mapsto|x-\mathbf{p}_j|\) sends \(C_j\) to an interval containing
\(0\) and \(|\mathbf{p}_j|>2\delta\). Thus every circle
\(\partial B_t(\mathbf{p}_j)\), \(\delta\le t\le2\delta\), meets
the zero set of \(u_j\).

Set \(K=\overline{B_{3\delta}(z)}\). Phase convergence gives
\(\lvert K\cap\{u_j=0\}\rvert\to0\). No circle in the
preceding range can be entirely zero for large \(j\):
subharmonicity and nonnegativity would imply \(u_j=0\) in
\(B_t(\mathbf{p}_j)\subset K\), whose area is at least
\(\pi\delta^2\). Each such circle therefore contains both a
zero point and a positive point, and hence meets \(\Gamma_j\).
The map \(x\mapsto|x-\mathbf{p}_j|\), which is Lipschitz with constant \(1\), consequently
sends \(\Gamma_j\cap K\) onto a set containing
\([\delta,2\delta]\). It follows that
\[
	\mathcal H^1(\Gamma_j\cap K)\ge\delta.
\]
The singular set has zero \(\mathcal H^1\) measure, also after
rescaling, so
\[
	\mathcal H^1\!\left(
	K\cap\partial_{\mathrm{reg}}\{u_j>0\}\right)\ge\delta.
\]

Integration by parts on each regular curve patch and the
one sided Bernoulli trace give
\[
	\Delta u_j\ge\sqrt{x_2}\,
	\mathcal H^1\mres\partial_{\mathrm{reg}}\{u_j>0\}
	\quad\hbox{in }\{x_2>0\}.
\]
Choose \(0\le\zeta\in C_c^\infty(B_1\cap\{x_2>0\})\)
equal to one on \(K\). Since \(\Delta u_j\ge0\), we infer
\[
	0<\delta\sqrt{\min_Kx_2}
	\le\int\zeta\,\dx{\Delta u_j}
	=\int u_j\Delta\zeta\,\dx{x}\longrightarrow0,
\]
contradicting the local uniform convergence \(u_j\to0\).
Therefore \(\Sigma^u=\varnothing\).
\end{proof}

\begin{proof}[Completion of Theorem~\ref{thm:stokes}\textup{(ii)}]
By Proposition~\ref{prop:no-full-density} and part~\textup{(i)},
every stagnation point has the Stokes blow-up, with convergence
of the functions and phases. It remains to show that
\(\mathcal S^u\) is locally finite.

Suppose that distinct stagnation points \(x^j\) converge to
\(x^0\in\Omega\). The set \(\mathcal S^u\) is relatively
closed, so \(x^0\in\mathcal S^u\). Set
\[
	r_j=|x^j-x^0|,\qquad
	\mathbf{q}_j=(x^j-x^0)/r_j\in\{(-1,0),(1,0)\}.
\]
After passing to a subsequence, \(\mathbf{q}_j=\mathbf{q}\) is fixed. Choose
\(s>0\) so small that \(\overline{B_{2s}(\mathbf{q})}\) is disjoint
from the closed Stokes sector. Scaling and \eqref{eq:boundary-Weiss-derivative-conclusion} give, for all sufficiently large \(j\),
\[
	\frac1{\sqrt3}
	=\Phi^{\mathrm{bound}}_{x^j,u}(0+)
	\le\Phi^{\mathrm{bound}}_{x^j,u}(r_js)
	=\Phi^{\mathrm{bound}}_{\mathbf{q},u_{r_j}}(s).
\]
The boundary term in the Weiss functional is nonpositive.
Consequently,
\[
	\Phi^{\mathrm{bound}}_{\mathbf{q},u_{r_j}}(s)
	\le s^{-3}\int_{B_s(\mathbf{q})}
	\bigl(|\nabla u_{r_j}|^2
	+x_2\chi_{\{u_{r_j}>0\}}\bigr)\dx{x}
	\longrightarrow0.
\]
The limit follows from strong gradient convergence and phase
convergence on \(B_{2s}(\mathbf{q})\), where \(\US=0\). This contradicts
the positive lower bound and proves local finiteness.
\end{proof}

We now impose the curve and trace conditions of
Theorem~\ref{thm:stokes}\textup{(iii)}. The next lemma recovers the
coefficient one asymptotically. The mean-value inequality for subharmonic
functions controls \((|\nabla u_r|^2-x_2)^+\), and a harmonic barrier gives an error
proportional to height. Together with the flatness argument below,
this provides uniform estimates near stagnation without bounds on
the initial curve norms.

\begin{lemma}\label{lem:pressure-excess}
Let \(u\) be a variational solution of \eqref{eq.problem} with
\(\Gamma=\partial_{\mathrm{reg}}\{u>0\}\), as defined in
Remark~\ref{rem:one-sided-trace}. Fix \(x^0\in\mathcal S^u\), and
suppose that \(u_r(x)=r^{-3/2}u(x^0+rx)\) converges strongly to
\(\US\) in \(W^{1,2}_{\loc}(\R^2)\) as \(r\to0^+\).
For every \(L>0\), there are \(\varepsilon_r=\varepsilon(r;L)\ge0\),
tending to zero, such that
\begin{equation}\label{eq:relative-pressure-conclusion}
	|\nabla u_r|^2\le(1+\varepsilon_r)x_2
	\quad\hbox{in }B_L\cap\{u_r>0\}.
\end{equation}
\end{lemma}

\begin{proof}
Fix \(L>0\) and take \(r\) small enough that \(u_r\) is defined
on \(B_{4L}\). By \eqref{eq:main-gradient}, we obtain
\(|\nabla u_r|^2\le Cx_2\) in the positive phase, uniformly in
\(r\). Define on \(B_{4L}\)
\[
	w(x)=
	\begin{cases}
		(|\nabla u_r(x)|^2-x_2)^+,&u_r(x)>0,\\
		0,&u_r(x)=0.
	\end{cases}
\]
This function is continuous. At a free boundary point above the
line \(\{x_2=0\}\), this follows from the one sided \(C^1\)
trace and the Bernoulli condition. At the line, it follows from
\(0\le w(x)\le Cx_2^+\). Harmonicity gives
\begin{align*}
	\Delta\bigl(|\nabla u_r|^2-x_2\bigr)
	&=2\sum_{i,j=1}^2(\partial_{x_i x_j}u_r)^2
	+2\sum_{i=1}^2\partial_{x_i}u_r\,
	\partial_{x_i}\Delta u_r-\Delta x_2\\
	&=2\sum_{i,j=1}^2(\partial_{x_i x_j}u_r)^2\ge0.
\end{align*}
Hence \(w\) is subharmonic in \(\{u_r>0\}\).

To verify subharmonicity across the interface, take any disk
\(D\Subset B_{4L}\) and let \(H\) be harmonic in \(D\),
continuous on \(\overline D\), and equal to \(w\) on
\(\partial D\). Since \(w\ge0\), the maximum principle gives
\(H\ge0\). On each connected component of
\(D\cap\{u_r>0\}\), the function \(w-H\) is subharmonic.
At boundary points in \(D\), we have \(w=0\le H\); at
boundary points on \(\partial D\), we have \(w=H\).
The maximum principle therefore gives \(w\le H\) in every
component. On the zero phase, the same inequality follows from
\(w=0\) and \(H\ge0\). The harmonic comparison criterion now
shows that \(w\) is subharmonic throughout \(B_{4L}\).

Equation \eqref{eq:stokes-profile} gives \(|\nabla\US|^2=\rho/2\le x_2\)
in its positive sector, and the gradient vanishes almost everywhere
outside it. Hence, almost everywhere in \(\{u_r>0\}\), we have
\[
	w\le\bigl||\nabla u_r|^2-|\nabla\US|^2\bigr|
	\le (|\nabla u_r|+|\nabla\US|)
	|\nabla u_r-\nabla\US|.
\]
The Cauchy--Schwarz inequality and strong convergence give
\[
	\int_{B_{4L}}w\,\dx{x}
	\le\bigl(\|\nabla u_r\|_{L^2(B_{4L})}
	+\|\nabla\US\|_{L^2(B_{4L})}\bigr)
	\|\nabla u_r-\nabla\US\|_{L^2(B_{4L})}
	\longrightarrow0.
\]
For \(x\in B_{3L}\), the mean-value inequality for subharmonic functions on
\(B_{L/2}(x)\subset B_{4L}\) yields
\[
	S:=\sup_{B_{3L}}w
	\le\frac{4}{\pi L^2}\int_{B_{4L}}w\,\dx{x}
	\longrightarrow0.
\]

Let \(h\) be bounded and harmonic in \(B_{2L}\cap\{x_2>0\}\),
with boundary values one on the open semicircle and zero on the open
diameter. The maximum principle gives \(0\le h\le1\). Odd reflection
and the interior gradient estimate give \(|\nabla h|\le C/L\)
in \(B_{3L/2}\), so
\[
	0\le h(x_1,x_2)
	=\int_0^{x_2}\partial_{x_2}h(x_1,t)\,\dx{t}
	\le Cx_2/L
	\quad\hbox{for }x\in B_L\cap\{x_2>0\}.
\]
The function \(w-Sh\) is subharmonic in the upper half disk.
On the open semicircle we have \(w\le S\) and \(h\to1\); on the
open diameter, \(w\le Cx_2\to0\) and \(h\to0\); and at its
two endpoints, \(w\to0\), while \(h\ge0\). Thus, at every
point \(\mathbf{p}\) of the boundary, we have
\[
	\limsup_{\substack{x\to \mathbf{p}\\x\in B_{2L}\cap\{x_2>0\}}}
	\bigl(w(x)-Sh(x)\bigr)\le0.
\]
The maximum principle gives \(w\le Sh\) in the half disk and
therefore
\[
	(|\nabla u_r(x)|^2-x_2)^+
	\le \frac{CS}{L}x_2
	\quad\hbox{for }x\in B_L\cap\{u_r>0\}.
\]
Taking \(\varepsilon_r=CS/L\) proves
\eqref{eq:relative-pressure-conclusion}, since \(S\to0\).
The constants may depend on \(L\).
\end{proof}

For the interface estimates, fix \(x^0\in\mathcal S^u\), an open
set \(U\subset\Omega\), and \(\rho>0\) such that
\begin{equation}\label{eq:localization-geometry}
	\overline{B_{16\rho}(x^0)}\subset U,\qquad
	B_{16\rho}(x^0)\cap\{u=0\}\subset K_0,
\end{equation}
where \(K_0\) is the connected component of \(U\cap\{u=0\}\)
containing \(x^0\). Such a choice is possible under the hypotheses
of Theorem~\ref{thm:stokes}\textup{(iii)}: choose \(U\Subset\Omega\)
with finitely many components of the zero phase, which are relatively closed,
and a disk avoiding those not containing \(x^0\). By part~\textup{(ii)},
we may also arrange that \(x^0\) is its only stagnation point.
Decrease \(\rho\) further so that
\(B_{256\rho}(x^0)\Subset\Omega\).\footnote{By \eqref{eq:main-gradient},
Theorem~\ref{thm:growth} controls \(\nabla u\) only on \(B_{R/16}(x^0)\).
Choosing \(R=256\rho\) gives \(R/16=16\rho\), the radius needed in
\eqref{eq:localization-growth}; any \(R\ge256\rho\) would serve equally well.}
Applying Theorem~\ref{thm:growth} on this larger ball gives \(M<\infty\)
such that
\begin{equation}\label{eq:localization-growth}
	|\nabla u|^2\le Mx_2^+
	\quad\hbox{a.e. in }B_{16\rho}(x^0).
\end{equation}
The next two estimates are uniform over free boundary points in
\(B_\rho(x^0)\).

\begin{lemma}\label{lem:arc-crossing}
Let \(u\) be a variational solution of \eqref{eq.problem} with
\(\Gamma=\partial_{\mathrm{reg}}\{u>0\}\). Fix
\(x^0\in\mathcal S^u\), an open set \(U\subset\Omega\),
and \(\rho>0\) satisfying \eqref{eq:localization-geometry}.
If \(\mathbf{b}=(b_1,h)\in\partial\{u>0\}\cap B_\rho(x^0)\), \(h>0\),
and \(0<s\le h/2\), then the free boundary curve through \(\mathbf{b}\)
contains an arc reaching \(\partial B_s(\mathbf{b})\), and
\begin{equation}\label{eq:arc-mass-conclusion}
	\Delta u(B_s(\mathbf{b}))\ge s\sqrt{h-s}
	\ge 2^{-1/2}s\sqrt h.
\end{equation}
\end{lemma}

\begin{proof}
The disk \(\overline{B_s(\mathbf{b})}\) lies in
\(B_{2\rho}(x^0)\cap\{x_2>0\}\), where every free boundary point has
a single regular curve patch, i.e.\ there is a neighborhood in
which the interface is an embedded \(C^{2,\alpha}\) curve with the positive
phase on one side (see Remark~\ref{rem:one-sided-trace}).
Continue the curve through \(\mathbf{b}\).
The interface is relatively closed, so every accumulation point in
this closed disk belongs to it. The single patch property places
such a point on the same component and permits continuation; it
excludes an interior endpoint or a distinct accumulating sheet.
At a boundary accumulation point the curve therefore reaches
\(\partial B_s(\mathbf{b})\). Otherwise the component is closed
in the disk and hence is a compact embedded one-dimensional
manifold without boundary, so it is a Jordan curve.

Let \(D\subset B_s(\mathbf{b})\) be its bounded interior, so
\(x^0\notin\overline D\). The positive phase lies on one consistent
side. On the interior side, a positive harmonic component would have
compact closure in \(\Omega\) and zero boundary values, contradicting
the maximum principle. On the exterior side, a positive collar
surrounds \(\partial D\). By \eqref{eq:localization-geometry}, we deduce
\(\partial D\subset K_0\). Thus \(K_0\cap\overline D\) is
nonempty and relatively closed in \(K_0\). It is also relatively
open: this is immediate at points of \(D\), while the positive
exterior collar excludes zero points outside \(\overline D\)
near \(\partial D\). Connectedness of \(K_0\) therefore gives
\(K_0\subset\overline D\), contradicting
\(x^0\in K_0\setminus\overline D\). Thus the curve reaches
\(\partial B_s(\mathbf{b})\).

The arc from \(\mathbf{b}\) to its first exit has length at least \(s\).
Integration by parts and the Bernoulli condition give
\begin{equation}\label{eq:regular-boundary-measure}
	\Delta u=\sqrt{x_2}\,
	\mathcal H^1\mres\partial\{u>0\}
	\quad\text{in }B_s(\mathbf{b});
\end{equation}
see also \cite[Lemma~3.4]{VW11}.
Integrating over the arc and using \(x_2\ge h-s\) gives
\eqref{eq:arc-mass-conclusion}.
\end{proof}

The mass bound \eqref{eq:arc-mass-conclusion} also gives a lower bound
for the Weiss functional at the vertical projection of each interface point.

\Needspace{8\baselineskip}
\begin{lemma}\label{lem:Weiss-detection}
Let \(u,x^0,U,\rho\) satisfy the assumptions of
Lemma~\ref{lem:arc-crossing}, and let \(M>0\) be a constant such that
\[
	|\nabla u|^2\le Mx_2^+
	\quad\text{a.e. in }B_{16\rho}(x^0).
\]
There is \(\beta=\beta(u,x^0,U,\rho,M)>0\) such that, for every
\(\mathbf{b}=(b_1,h)\in\partial\{u>0\}\cap B_\rho(x^0)\) with \(h>0\),
the projection \(\mathbf{p}=(b_1,0)\) satisfies
\begin{equation}\label{eq:Weiss-detection-conclusion}
	\Phi^{\mathrm{bound}}_{\mathbf{p},u}(4h)\ge\beta.
\end{equation}
\end{lemma}

\begin{proof}
The inequalities \(h\le\rho\) and \(|\mathbf{p}-x^0|\le\rho\) give
\begin{equation}\label{eq:detection-disk-inclusions}
	B_{8h}(\mathbf{p})\subset B_{9\rho}(x^0)
	\subset B_{16\rho}(x^0).
\end{equation}
If the conclusion fails for fixed \(u,x^0,U,\rho,M\), choose
\(\mathbf{b}_j=(b_{j,1},h_j)\) such that its projection
\(\mathbf{p}_j=(b_{j,1},0)\) satisfies
\(\Phi^{\mathrm{bound}}_{\mathbf{p}_j,u}(4h_j)\to0\). Set
\[
	v_j(y)=h_j^{-3/2}u(\mathbf{p}_j+h_jy),\qquad y\in B_8.
\]
These functions are defined on the common disk \(B_8\) by
\eqref{eq:detection-disk-inclusions}, regardless of the limit of
\(h_j\). The bound \eqref{eq:localization-growth} and vertical
integration from the line \(\{x_2=0\}\) give
\[
	|\nabla v_j|^2\le My_2^+,\qquad
	0\le v_j(y)\le\frac23\sqrt M\,(y_2^+)^{3/2}.
\]
A subsequence converges locally uniformly and weakly in
\(W^{1,2}_{\loc}(B_8)\) to \(v\), which vanishes below the line
and is harmonic in its positive phase. Truncation gives
\eqref{eq:energy-conclusion} for \(v\); comparison with the identities
for \(v_j\) gives strong \(W^{1,2}_{\loc}(B_8)\) convergence.

The phase indicators also have uniform local \(BV\) bounds above
the line. Indeed, testing \(\Delta v_j\ge0\) against a cutoff function
supported in \(B_7\) and equal to one on \(B_6\) gives, uniformly
in the moving centers,
\[
	\Delta v_j(B_6)\le C\sqrt M.
\]
By \eqref{eq:regular-boundary-measure} on the regular curve patches, we obtain
\[
	\int_{B_6\cap\{y_2>0\}}\sqrt{y_2}\,
	\dx{|D\chi_{\{v_j>0\}}|}\le C\sqrt M.
\]
The weight is bounded below on \(\{y_2>\tau\}\) for each
\(\tau>0\). By \(BV\) compactness, we extract a diagonal subsequence as
\(\tau\to0^+\); the vanishing area of the remaining strip then gives
\[
	\chi_{\{v_j>0\}}\longrightarrow\chi_0
	\quad\hbox{in }L^1_{\loc}(B_6),\qquad
	\chi_0\ge\chi_{\{v>0\}}.
\]
Here \(\chi_0\) is a characteristic function, zero below the line;
its value on \(\{v=0\}\) remains undetermined.
Lemma~\ref{lem:arc-crossing} and scaling give, with \(\mathbf{e}_2=(0,1)\),
\[
	\Delta v_j(B_{1/2}(\mathbf{e}_2))
	=h_j^{-3/2}\Delta u(B_{h_j/2}(\mathbf{b}_j))\ge 2^{-3/2}.
\]
Take \(0\le\eta\in C_c^\infty(B_2)\) equal to one on
\(\overline{B_{1/2}(\mathbf{e}_2)}\). Passing to the limit in distributions gives
\begin{equation}\label{eq:detection-nonzero-limit}
	\int_{B_2}v\Delta\eta\dx{y}
	=\lim_{j\to\infty}\int\eta\,\dx{\Delta v_j}
	\ge 2^{-3/2}.
\end{equation}
In particular, \(v\not\equiv0\) in \(B_2\).
The boundary Weiss functional is nonnegative at every center on the line \(\{x_2=0\}\): use \cite[Lemma~4.2(ii)]{VW11} and
\eqref{eq:boundary-Weiss-derivative-conclusion} at free boundary points;
at other points \(u\) vanishes nearby, so its limiting value is zero.
Set
\[
	W_j(s)=\Phi^{\mathrm{bound}}_{0,v_j}(s)
	=\Phi^{\mathrm{bound}}_{\mathbf{p}_j,u}(h_js).
\]
By the choice of \(\mathbf{p}_j,h_j\) and \eqref{eq:boundary-Weiss-derivative-conclusion}, we obtain
\begin{equation}\label{eq:detection-Weiss-squeeze}
	0\le W_j(s)\le W_j(4)\longrightarrow0
	\qquad(0<s<4).
\end{equation}
From the strong convergence of the gradient, the local uniform convergence,
and the convergence of the phase indicators, we obtain, for almost every
\(s\in(0,4)\),
\[
	W_j(s)\longrightarrow W_{v,\chi_0}(s),
\]
where \(W_{v,\chi_0}\) is the functional
\eqref{eq:boundary-Weiss-conclusion} with phase term \(y_2\chi_0\).
All spheres lie inside the common domain, so interior convergence
suffices. Thus \(W_{v,\chi_0}(s)=0\) for almost every \(s\in(0,4)\).
By the coarea formula, \eqref{eq:boundary-Weiss-integrated-conclusion} reads
\[
	W_j(t)-W_j(s)
	=2\int_{B_t\setminus B_s}|y|^{-5}
	\left(y\cdot\nabla v_j-\frac32v_j\right)^2\dx{y},
\]
for \(0<s<t<4\). By \eqref{eq:detection-Weiss-squeeze} and strong
convergence on annuli, we infer that \(y\cdot\nabla v=3v/2\) almost everywhere in
\(B_4\). This homogeneity and \eqref{eq:energy-conclusion} give,
for almost every \(s<4\),
\[
	\int_{B_s}|\nabla v|^2\dx{y}
	=\frac{3}{2s}\int_{\partial B_s}v^2\dH{1},\qquad
	0=W_{v,\chi_0}(s)
	=s^{-3}\int_{B_s}y_2\chi_0\dx{y}.
\]
Since \(\chi_0\ge\chi_{\{v>0\}}\) and both functions vanish below
the line, the last equality forces \(v=0\) in \(B_4\), contradicting
\eqref{eq:detection-nonzero-limit}. The contradiction fixes a positive
lower bound for the given \(u,x^0,U,\rho,M\), uniformly in \(\mathbf{b},h\).
\end{proof}

We can now locate the free boundary on annuli about a stagnation
point. Write
\[
	u_r(x)=r^{-3/2}u(x^0+rx),\qquad
	\Gamma_r=\partial\{u_r>0\},\qquad
	\Gamma_{\mathrm{st}}=\partial\{\US>0\},
\]
and set \(\mathcal A_{a,A}=\{a\le|x|\le A\}\).

\Needspace{10\baselineskip}
\begin{lemma}
Let \(u\) satisfy the hypotheses of
Theorem~\ref{thm:stokes}\textup{(iii)}, and fix
\(x^0\in\mathcal S^u\). For every \(0<a<A<\infty\), the following limits hold:
\begin{equation}\label{eq:interface-Hausdorff-conclusion}
	\begin{split}
		&\sup_{\mathbf{b}\in\Gamma_r\cap\mathcal A_{a,A}}
		\operatorname{dist}(\mathbf{b},\Gamma_{\mathrm{st}})\longrightarrow0,\\
		&\sup_{\mathbf{b}\in\Gamma_{\mathrm{st}}\cap\mathcal A_{a,A}}
		\operatorname{dist}(\mathbf{b},\Gamma_r)\longrightarrow0
		\qquad\hbox{as }r\to0^+.
	\end{split}
\end{equation}
\end{lemma}

The distances in \eqref{eq:interface-Hausdorff-conclusion} are to
the full interfaces, regardless of their intersections with the
edges of the annulus.

\begin{proof}
Choose \(U,\rho,M\) satisfying \eqref{eq:localization-geometry}
and \eqref{eq:localization-growth} so that \(x^0\) is the only stagnation
point in \(B_{16\rho}(x^0)\). Fix \(\beta>0\) from
Lemma~\ref{lem:Weiss-detection}. For the annulus \(\mathcal A_{a,A}\),
take \(0<r<\rho/(A+1)\), so that
\(x^0+r\mathcal A_{a,A}\subset B_\rho(x^0)\).

By Theorem~\ref{thm:stokes}\textup{(ii)}, \(u_r\) converges to \(\US\) locally
uniformly and strongly in \(W^{1,2}_{\loc}\), and the phase indicators
converge in \(L^1_{\loc}\). If the first limit in
\eqref{eq:interface-Hausdorff-conclusion} fails, choose
\(\mathbf{b}_j\in\Gamma_{r_j}\cap\mathcal A_{a,A}\), \(r_j\to0^+\),
with \(\mathbf{b}_j\to \mathbf{b}\notin\Gamma_{\mathrm{st}}\). Then \(b_2\ge0\).
Local uniform convergence excludes \(\US(\mathbf{b})>0\), since it would
make \(u_{r_j}\) positive near \(\mathbf{b}\).

If \(b_2>0\), choose \(0<d<b_2/8\) such that \(B_{4d}(\mathbf{b})\)
lies above the line and inside \(\{\US=0\}\).
Applying Lemma~\ref{lem:arc-crossing} at
\(x^0+r_jb_j\in B_\rho(x^0)\) with radius \(r_jd\) gives
\[
	\Delta u_{r_j}(B_d(\mathbf{b}_j))
	\ge d\sqrt{(\mathbf{b}_j)_2-d}\ge c(\mathbf{b},d)>0.
\]
Testing with a nonnegative cutoff function supported in \(B_{3d}(\mathbf{b})\) and
equal to one on \(B_{2d}(\mathbf{b})\) contradicts
\(\Delta u_{r_j}\rightharpoonup\Delta\US=0\) there.

If \(b_2=0\), then \(|\mathbf{b}|\ge a>0\) and isolation of \(x^0\)
gives \(t_j=(\mathbf{b}_j)_2>0\) for large \(j\), with \(t_j\to0\).
Set \(\mathbf{q}_j=((\mathbf{b}_j)_1,0)\to \mathbf{q}=\mathbf{b}\). Applying
Lemma~\ref{lem:Weiss-detection}
(cf.\ equation~\eqref{eq:Weiss-detection-conclusion}) to the fixed
solution \(u\) at \(x^0+r_jb_j\) and rescaling gives
\[
	\Phi^{\mathrm{bound}}_{\mathbf{q}_j,u_{r_j}}(4t_j)\ge\beta.
\]
Here \(\beta\) is independent of \(j\), since the original local
data are fixed. Choose \(s>0\) with \(\overline{B_{2s}(\mathbf{q})}\)
disjoint from the closed Stokes sector. For large \(j\),
\(4t_j<s\) and \(B_s(\mathbf{q}_j)\subset B_{2s}(\mathbf{q})\) hold. By
\eqref{eq:boundary-Weiss-derivative-conclusion} and
\eqref{eq:boundary-Weiss-conclusion}, we obtain
\[
	\begin{split}
		\beta
		&\le\Phi^{\mathrm{bound}}_{\mathbf{q}_j,u_{r_j}}(s)\\
		&\le s^{-3}\int_{B_s(\mathbf{q}_j)}
		\bigl(|\nabla u_{r_j}|^2
		+x_2\chi_{\{u_{r_j}>0\}}\bigr)\dx{x}
		\longrightarrow0.
	\end{split}
\]
The convergence follows from the strong gradient and phase convergence
on \(B_{2s}(\mathbf{q})\), where \(\US=0\). This contradiction proves the
first limit in \eqref{eq:interface-Hausdorff-conclusion}.

For the reverse inclusion, every sufficiently small disk centered on
\(\Gamma_{\mathrm{st}}\setminus\{0\}\) contains a point inside the
Stokes sector and a smaller disk in the adjacent air. Local uniform
and phase convergence supply, respectively, a positive point and a
zero of \(u_r\) in the original disk. The segment joining them meets
\(\Gamma_r\). A finite covering of
\(\Gamma_{\mathrm{st}}\cap\mathcal A_{a,A}\) gives the second limit.
\end{proof}

For convergence of the normals, we use the following form of
De~Silva's flatness theorem \cite[Theorem~1.1]{DS11}.
Appendix~\ref{app:flatness} verifies the viscosity hypotheses
and tracks the Bernoulli coefficient through the flatness
iteration to obtain the stated estimate for the normal.

\begin{lemma}\label{lem:local-water-wave-flatness}
Let \(u\in C^0(B_1)\) be nonnegative and harmonic in \(\{u>0\}\),
with \(0\in\partial\{u>0\}\). At every point of
\(B_1\cap\partial\{u>0\}\), assume that the interface is a
\(C^2\) curve patch with the positive phase on one side, and that
\(u\) has a \(C^1\) trace from that side satisfying
\[
	|\nabla u|=Q\quad\hbox{on }B_1\cap\partial\{u>0\},
	\qquad Q\in C^{0,1}(B_1),\quad Q>0,\quad Q(0)=1.
\]
There are universal constants \(\varepsilon_*,C_*>0\) and
\(\alpha_*\in(0,1)\) such that, if \(n\) is a unit vector, \(0<\varepsilon\le\varepsilon_*\), and
\begin{equation}\label{eq:local-flatness-sandwich}
	(z\cdot n-\varepsilon)^+\le u(z)
	\le(z\cdot n+\varepsilon)^+\quad(z\in B_1),
	\qquad [Q]_{C^{0,1}(B_1)}\le\varepsilon^2,
\end{equation}
then the interface is a \(C^{1,\alpha_*}\) graph in \(B_{1/2}\),
and its outward unit normal \(\nu\) satisfies
\[
	|\nu(0)+n|\le C_*\varepsilon.
\]
\end{lemma}

\begin{proof}[Proof of Theorem~\ref{thm:stokes}~\textup{(iii)}]
Let \(u\) satisfy the stated assumptions. Fix \(x^0\in\mathcal S^u\)
and write \(u_r(x)=r^{-3/2}u(x^0+rx)\). By part~\textup{(ii)}, we deduce that
\(u_r\to\US\) locally uniformly, with convergence of the phases.
We use \eqref{eq:interface-Hausdorff-conclusion} and
\eqref{eq:relative-pressure-conclusion} to apply
Lemma~\ref{lem:local-water-wave-flatness} uniformly on a fixed annulus.

Set
\[
	A=\{1/2\le|x|\le2\},\qquad
	A'=\{1/4\le|x|\le3\},\qquad
	\Gamma_r=\partial\{u_r>0\}.
\]
For small \(r\), each \(\mathbf{b}\in\Gamma_r\cap A\) lies near one Stokes
ray, away from the other, and \(b_2\ge m:=1/8\). Let \(\eta\)
be the outward normal to that ray and \(n=-\eta\), so
\[
	\eta\in\{(-1/2,-\sqrt3/2),(1/2,-\sqrt3/2)\}.
\]
Fix \(0<\varepsilon\le\min\{\varepsilon_*,1/4\}\). Choose
\(\delta>0\) so that each \(B_{2\delta}(\mathbf{b})\) lies in
\(\operatorname{int}A'\), meets only the neighborhood of its
corresponding ray, and satisfies
\begin{equation}\label{eq:flatness-radius-choice}
	\delta\le1/128,\qquad \frac{\delta}{m}\le\frac{\varepsilon^2}{64}.
\end{equation}
We may decrease \(\delta\) once more, still subject to
\eqref{eq:flatness-radius-choice}, before choosing \(r\).
Let \(\bar{\mathbf{b}}\) be the orthogonal projection of \(\mathbf{b}\) onto the ray;
for small \(r\), \(|\mathbf{b}-\bar{\mathbf{b}}|\le\delta\) uniformly in \(\mathbf{b}\).
On the positive side of the fixed annulus, \(\US\) has bounded
second derivatives, zero boundary value, and inward derivative
\(\sqrt{\bar b_2}\). Taylor expansion and positive parts give
\[
	\left|\US(x)-\sqrt{\bar b_2}
	((x-\bar{\mathbf{b}})\cdot n)^+\right|
	\le C|x-\bar{\mathbf{b}}|^2
	\quad\hbox{in }B_{2\delta}(\mathbf{b}).
\]
Here and below \(C\) is uniform in \(\mathbf{b}\); both terms vanish on
the air side. Since the positive part map is Lipschitz and
\(|\sqrt{b_2}-\sqrt{\bar b_2}|\le C|\mathbf{b}-\bar{\mathbf{b}}|\), we obtain
\begin{equation}\label{eq:stokes-patch-Taylor}
	\sup_{B_{2\delta}(\mathbf{b})}
	\left|\US(x)-\sqrt{b_2}((x-\mathbf{b})\cdot n)^+\right|
	\le C\delta^2+C|\mathbf{b}-\bar{\mathbf{b}}|.
\end{equation}
Decrease \(\delta\) so that
\(C\delta^2/(\delta\sqrt m)\le\varepsilon/16\), then fix it.
By \eqref{eq:interface-Hausdorff-conclusion} and uniform convergence
on \(A'\), choose \(r\) so that the sum of the second error in
\eqref{eq:stokes-patch-Taylor} and
\(\|u_r-\US\|_{L^\infty(A')}\) is at most
\(\delta\sqrt m\,\varepsilon/16\), uniformly in \(\mathbf{b}\).
The functions
\[
	w_{r,\mathbf{b}}(z)=\frac{u_r(\mathbf{b}+\delta z)}{\delta\sqrt{b_2}},
	\qquad z\in B_2,
\]
satisfy, uniformly in \(\mathbf{b}\in\Gamma_r\cap A\),
\begin{equation}\label{eq:patch-plane-approximation}
	\left|w_{r,\mathbf{b}}(z)-(z\cdot n)^+\right|
	\le\varepsilon/8\qquad(z\in B_2).
\end{equation}
Nonnegativity and \eqref{eq:patch-plane-approximation} give the lower
bound in \eqref{eq:local-flatness-sandwich}.
For the upper bound, decrease \(r\) so that
\eqref{eq:interface-Hausdorff-conclusion} places every interface
point of \(B_{2\delta}(\mathbf{b})\) in
\[
	|(x-\mathbf{b})\cdot n|<\delta\varepsilon/16.
\]
The connected cap \(B_{3/2}\cap\{z\cdot n\le-\varepsilon/8\}\)
then contains no interface point. Require also the Hausdorff error
on \(A'\) to be less than \(\delta\varepsilon/64\). Its image
lies in the limiting air, since for \(x=\mathbf{b}+\delta z\) in that image,
\[
	(x-\bar{\mathbf{b}})\cdot n
	\le |\mathbf{b}-\bar{\mathbf{b}}|-\delta\varepsilon/8
	<-7\delta\varepsilon/64.
\]
The image contains \(B_{\delta/4}(\mathbf{b}-\delta n)\), whose area is
\(\pi\delta^2/16\), uniformly in \(\mathbf{b}\). Require the
\(L^1(A')\) error of the phase indicators to be smaller than this
area. Phase convergence supplies a zero in every cap;
connectedness then makes each cap a zero region. Thus, uniformly
in all centers, we obtain
\begin{equation}\label{eq:local-flatness-air}
	w_{r,\mathbf{b}}=0\quad\hbox{in }
	B_{3/2}\cap\{z\cdot n\le-\varepsilon/8\}.
\end{equation}
By \eqref{eq:relative-pressure-conclusion} on a fixed disk containing
\(A'\), we may also require \(\varepsilon_r\le\varepsilon^2/64\).
With \(\kappa_{\mathbf{b}}=\delta/b_2\), we estimate
\[
	|\nabla w_{r,\mathbf{b}}|^2
	\le(1+\varepsilon_r)(1+\kappa_{\mathbf{b}} z_2)
	\le(1+\varepsilon^2/64)(1+\varepsilon^2/32)
	\le(1+\varepsilon/8)^2\quad\hbox{a.e. in }B_2.
\]
The last inequality uses \(\varepsilon\le1/4\); the gradient
vanishes almost everywhere in the zero phase. Convexity of \(B_2\)
therefore makes \(w_{r,\mathbf{b}}\) Lipschitz with constant \(1+\varepsilon/8\) there.

For \(z\in B_1\) with \(t=z\cdot n\ge-\varepsilon/8\),
\eqref{eq:local-flatness-air} gives a zero at
\(z-(t+\varepsilon/8)n\in B_{3/2}\). The connecting segment stays
in \(B_{3/2}\), so
\[
	w_{r,\mathbf{b}}(z)\le(1+\varepsilon/8)(t+\varepsilon/8)
	\le t+\varepsilon.
\]
For \(t< -\varepsilon/8\), \(w_{r,\mathbf{b}}=0\), proving the upper
bound in \eqref{eq:local-flatness-sandwich} on \(B_1\).
The Bernoulli coefficient is
\[
	Q_{\mathbf{b}}(z)=\sqrt{1+\kappa_{\mathbf{b}} z_2},\qquad Q_{\mathbf{b}}(0)=1.
\]
It is positive on \(B_1\), with
\begin{equation}\label{eq:normalized-water-wave-coefficient}
	[Q_{\mathbf{b}}]_{C^{0,1}(B_1)}
	\le\frac{\kappa_{\mathbf{b}}}{2\sqrt{1-\kappa_{\mathbf{b}}}}
	\le\kappa_{\mathbf{b}}\le\varepsilon^2/64.
\end{equation}
Since the patches lie above the line \(\{x_2=0\}\), the assumed curve
regularity and Bernoulli trace persist under rescaling. Thus
Lemma~\ref{lem:local-water-wave-flatness} applies at every center
with constants independent of \(r,\mathbf{b},\delta\); it gives
\[
	\sup_{\mathbf{b}\in\Gamma_r\cap A}|\nu_r(\mathbf{b})-\eta|\le C_*\varepsilon,
\]
where \(\eta\) is the normal to the nearest ray. Letting
\(\varepsilon\to0^+\) gives uniform convergence of the normals
on the annulus. Rescaling each interface point by its distance from
\(x^0\) gives the corresponding convergence as it approaches \(x^0\).

To count the branches, translate \(x^0\) to the origin and choose
\(\rho>0\) so small that, in
\[
	R_\rho=(-\rho,\rho)\times(0,\rho)
\]
every interface normal has vertical component negative and bounded
away from zero. Equation \eqref{eq:interface-Hausdorff-conclusion}
gives \(0<x_2\le2|x_1|\), excluding the vertical axis. Uniform
convergence of \(u_\rho\) makes the top edge positive, since
\(\{(t,1):|t|\le1\}\) lies inside the Stokes sector. Each bottom
point \((t,0)\), \(0<|t|<\rho\), has a neighborhood on which \(u=0\)
because the origin is an isolated stagnation point.

For each \(t\in(-\rho,0)\cup(0,\rho)\), the vertical segment
above \((t,0)\) meets the interface away from its endpoints.
The intersections are transverse and locally isolated, hence finite
by compactness. At every crossing the upward segment passes from
zero to positive values, so there is exactly one: a second would
require a passage back to the zero phase.

The interface is therefore one graph over each interval. This also
excludes closed loops and additional arcs, including those accumulating
at the origin, since they would give a further vertical crossing.
The height bound makes both graphs tend to the origin. Their
normals, and hence their slopes, tend to the Stokes values
\(-1/\sqrt3\) on the left and \(1/\sqrt3\) on the right.
Assigning to each graph the height zero at its endpoint, thereby
appending the origin, gives the claimed one sided \(C^1\) regularity.
\end{proof}

\subsubsection*{Acknowledgement}
\hyphenpenalty=10
\sloppy
L.~Du is supported by National Natural Science Foundation of China under Grants 12125102, 12526202, 12671245, 12631009. C.~Yang is supported by National Natural Science Foundation of China under Grant 12601435. F.~Zhang is supported by the Fundamental Research Program of Shanxi Province  under Grant 202303021222116.

\subsubsection*{Conflict of interest}
The authors declare that they have no conflict of interest.

\subsubsection*{AI declarations}
OpenAI's GPT was used to assist with language editing, citation checks, and the mathematical revision of the manuscript. The authors take full responsibility for the mathematical content and the final manuscript.

\subsubsection*{Data availability} No data was used in this research.

\FloatBarrier

\appendix

\section{Reduction to the approximation theorem}
\label{app:approximation-reduction}

We first state the approximation theorem used in the proof of
Theorem~\ref{thm:growth}, specialized to maps with values in
\(\C\simeq\R^2\). We use Theorem~1 of version~2 of
\cite{DHN26}. For a measurable function
\(W\) on \(B_1\), its associated boundary is
\[
	\mathcal F_W=\left\{x\in B_1:
	\begin{aligned}
		&|B_\rho(x)\cap\{W=0\}|\,|B_\rho(x)\cap\{W\ne0\}|>0\\
		&\hspace{35mm}\text{for every }0<\rho<1-|x|
	\end{aligned}\right\}.
\]
As in Section~\ref{sec:decay}, \((h)_B=|B|^{-1}\int_Bh\,\dx{x}\).

\begin{theorem}[De Philippis--Hirsch--Nahon, version~2 of \cite{DHN26}]
\label{thm:DHN-unit-disk}
Let \(K\ge0\) and \(W\in L^2(B_1;\C)\), with
\(\mathcal F_W\ne\varnothing\). Suppose that every disk
\(B_s(\mathbf{b})\subset B_1\) admits \(h_{B_s(\mathbf{b})}\in L^2(B_s(\mathbf{b});\C)\)
satisfying
\begin{align*}
	\int_{B_s(\mathbf{b})}|W-h_{B_s(\mathbf{b})}|^2\dx{x}
	&\le |B_s(\mathbf{b})\cap\{W=0\}|,\\
	s\|\nabla h_{B_s(\mathbf{b})}\|_{L^\infty(B_{s/2}(\mathbf{b}))}
	&\le K\left(\barint_{B_s(\mathbf{b})}
	|h_{B_s(\mathbf{b})}-(h_{B_s(\mathbf{b})})_{B_s(\mathbf{b})}|^2\dx{x}\right)^{1/2}.
\end{align*}
Then there are \(C,\alpha>0\), depending only on \(K\), such that
\[
	|W(x)|\le C\left(1+
	\operatorname{dist}(x,\mathcal F_W)^\alpha\|W\|_{L^2(B_1)}\right)
	\quad\text{for almost every }x\in B_{1/2}.
\]
\end{theorem}

We now derive the scaled form in Theorem~\ref{thm:DHN}.
The quoted theorem assumes \(\mathcal F_W\ne\varnothing\).
The additional assertion when \(\mathcal F_W=\varnothing\)
and the approximants are holomorphic is proved separately below.

Let \(K\ge0\), \(\Lambda>0\), and
\(W\in L^2(B_r(\mathbf{a});\C)\) satisfy
\eqref{eq:DHN-first}--\eqref{eq:DHN-second}. On the unit disk, set
\[
	\widetilde W(y)=\Lambda^{-1}W(\mathbf{a}+ry),\qquad
	\widetilde h_{B_t(\mathbf{c})}(y)
	=\Lambda^{-1}h_{B_{rt}(\mathbf{a}+r\mathbf{c})}(\mathbf{a}+ry).
\]
Here each disk \(B_t(\mathbf{c})\subset B_1\) corresponds to
\(B_{rt}(\mathbf{a}+r\mathbf{c})\subset B_r(\mathbf{a})\), and we identify \(\C\) with
\(\R^2\). Changing variables in \eqref{eq:DHN-first} gives
\begin{align*}
	\int_{B_t(\mathbf{c})}|\widetilde W-\widetilde h_{B_t(\mathbf{c})}|^2\dx{y}
	&=\frac1{r^2\Lambda^2}
	\int_{B_{rt}(\mathbf{a}+r\mathbf{c})}|W-h_{B_{rt}(\mathbf{a}+r\mathbf{c})}|^2\dx{x}\\
	&\le r^{-2}|B_{rt}(\mathbf{a}+r\mathbf{c})\cap\{W=0\}|\\
	&=|B_t(\mathbf{c})\cap\{\widetilde W=0\}|.
\end{align*}
Thus the error estimate has coefficient one after rescaling.
For the gradient estimate, use
\(\nabla_y\widetilde h=(r/\Lambda)\nabla_xh\) in
\eqref{eq:DHN-second} to obtain
\begin{align*}
	t\|\nabla\widetilde h_{B_t(\mathbf{c})}\|_{L^\infty(B_{t/2}(\mathbf{c}))}
	&=\frac{rt}{\Lambda}
	\|\nabla h_{B_{rt}(\mathbf{a}+r\mathbf{c})}\|_{L^\infty(B_{rt/2}(\mathbf{a}+r\mathbf{c}))}\\
	&\le K\left(\barint_{B_t(\mathbf{c})}
	|\widetilde h_{B_t(\mathbf{c})}
	-(\widetilde h_{B_t(\mathbf{c})})_{B_t(\mathbf{c})}|^2\dx{y}\right)^{1/2}.
\end{align*}
The two approximation hypotheses of Theorem~\ref{thm:DHN-unit-disk}
therefore hold on \(B_1\), with the same \(K\).

Assume first that \(\mathcal F_W\ne\varnothing\). The norm
and the associated boundary transform according to
\[
	\|\widetilde W\|_{L^2(B_1)}
	=\frac1{r\Lambda}\|W\|_{L^2(B_r(\mathbf{a}))},\qquad
	\mathcal F_{\widetilde W}=r^{-1}(\mathcal F_W-\mathbf{a}).
\]
Applying Theorem~\ref{thm:DHN-unit-disk} and multiplying its
conclusion by \(\Lambda\) gives the distance estimate in
Theorem~\ref{thm:DHN}. Since
\(\operatorname{dist}(x,\mathcal F_W)<3r/2\) for
\(x\in B_{r/2}(\mathbf{a})\), this also proves \eqref{eq:DHN-bound}.

Suppose next that \(\mathcal F_W=\varnothing\). By its
definition, every point has a disk on which either \(W=0\)
almost everywhere or \(\{W=0\}\) has measure zero. The points
of each type form an open set. These two open sets are disjoint:
if disks of opposite types intersected, their intersection would
have positive measure, giving a contradiction. They cover the
connected disk \(B_r(\mathbf{a})\), so only one type occurs. A countable
subcover then gives either \(W=0\) almost everywhere on
\(B_r(\mathbf{a})\), or \(|B_r(\mathbf{a})\cap\{W=0\}|=0\).

In the first case, \eqref{eq:DHN-bound} is immediate. In the
second, \eqref{eq:DHN-first} gives \(W=h_B\) almost everywhere
on every subdisk \(B\). If the approximants are holomorphic,
they agree on overlaps and define a holomorphic representative
of \(W\). Applying the mean-value inequality for subharmonic
functions to \(|W|^2\) yields
\[
	\|W\|_{L^\infty(B_{r/2}(\mathbf{a}))}
	\le Cr^{-1}\|W\|_{L^2(B_r(\mathbf{a}))}.
\]
This proves the remaining assertion of Theorem~\ref{thm:DHN}.
In its application to \(W=\Phi\),
Lemmas~\ref{lem:approximation} and~\ref{lem:decay} supply the
required integrability and approximation estimates, with
universal \(K\) and \(\Lambda=3r\).

\section{The local flatness implication}\label{app:flatness}

Lemma~\ref{lem:local-water-wave-flatness} follows from
De~Silva's flatness theorem \cite[Theorem~1.1]{DS11}. We verify that the classical boundary
conditions imply the viscosity conditions used there, and then
estimate the change of normal in the flatness iteration.

\begin{proof}[Proof of Lemma~\ref{lem:local-water-wave-flatness}]
Let \(u\) satisfy the assumptions of the lemma. The interior
equation holds classically. At a free boundary point \(\mathbf{p}\), let
\(n(\mathbf{p})\) denote the inward unit normal. The tangential derivative
of \(u\) vanishes, while its inward derivative equals \(Q(\mathbf{p})>0\).
The one sided \(C^1\) trace therefore gives
\[
	u(\mathbf{p}+y)=Q(\mathbf{p})(y\cdot n(\mathbf{p}))^++o(|y|).
\]
To check \cite[Definition~2.2]{DS11}, take \(\phi\in C^2\)
with \(\phi(\mathbf{p})=0\) and \(|\nabla\phi(\mathbf{p})|\ne0\). If
\(\phi^+\) touches \(u\) from below, inclusion of the positive
sides forces their inward directions to agree. Comparison along
\(n(\mathbf{p})\) then gives \(|\nabla\phi(\mathbf{p})|\le Q(\mathbf{p})\).
Contact from above gives the reverse inequality. Thus \(u\)
satisfies the required viscosity boundary conditions.

Rotate coordinates so that the direction \(n\) in
\eqref{eq:local-flatness-sandwich} is vertical. The equation
remains the Laplace equation, and the Bernoulli coefficient
satisfies \(Q(0)=1\) and
\([Q]_{C^{0,1/2}(B_1)}\le\sqrt2\,\varepsilon^2\).
After choosing a sufficiently small universal
\(\varepsilon_*\), \cite[Theorem~1.1]{DS11} implies that the
interface is a \(C^{1,\alpha_*}\) graph in \(B_{1/2}\).
It remains to compare its normal at the origin with \(-n\).

Fix a universal \(0<\theta\le1/4\) for which
\cite[Lemma~4.1]{DS11} applies, and decrease \(\varepsilon_*\)
to its threshold. The improvement of flatness, together with
the linear bound on the change of direction in its proof, gives
directions \(n_k\), starting with \(n_0=n\), such that
\begin{align*}
	(z\cdot n_k-\rho_k\varepsilon_k)^+
	&\le u(z)\le(z\cdot n_k+\rho_k\varepsilon_k)^+
	&&(z\in B_{\rho_k}),\\
	|n_{k+1}-n_k|&\le C\varepsilon_k,
\end{align*}
where \(\rho_k=\theta^k\) and
\(\varepsilon_k=2^{-k}\varepsilon\). To justify each successive
application, observe that the coefficient for
\(u(\rho_k z)/\rho_k\) satisfies
\[
	\|Q(\rho_k z)-1\|_{L^\infty(B_1)}
	\le\varepsilon^2\theta^k
	\le\varepsilon^2 4^{-k}=\varepsilon_k^2.
\]
Rotating to \(n_k\) preserves both this bound and the Laplace
equation. The coefficient is therefore small enough at every
stage of the iteration.

The changes of direction are summable, so \(n_k\) converges
to a unit vector \(n_\infty\), with
\[
	|n_\infty-n|\le C\sum_{k=0}^\infty\varepsilon_k
	\le2C\varepsilon,
	\qquad
	|n_\infty-n_k|\le2C\varepsilon_k.
\]
The slabs have width \(2\rho_k\varepsilon_k=o(\rho_k)\).
Consequently, their limiting direction is the inward normal
at the origin, and \(|\nu(0)+n|\le2C\varepsilon\).
The geometric decay allows a H\"older exponent no larger than
\(\log2/|\log\theta|\).

For the normalized water wave patches used in the proof of
Theorem~\ref{thm:stokes}\textup{(iii)}, the coefficient is
\(Q_{\mathbf{b}}(z)=\sqrt{1+\kappa_{\mathbf{b}} z_2}\). The choice of patch radius
makes \(\kappa_{\mathbf{b}}\) small, and
\eqref{eq:normalized-water-wave-coefficient} gives
\([Q_{\mathbf{b}}]_{C^{0,1}(B_1)}\le\kappa_{\mathbf{b}}\).
Thus the constants in the lemma are independent of the physical
height and the patch radius.
\end{proof}

\section{Use of the geometric classification}
\label{app:hypothesis-audit}

We explain how the assumptions in the three parts of
Theorem~\ref{thm:stokes} enter the proof. The common starting
point is Definition~\ref{def:variational}. Continuity,
\(C^2\) regularity in the positive phase, and the domain variation
identity \eqref{eq:stationarity}
imply harmonicity there and the energy identity
\eqref{eq:energy-conclusion}. Theorem~\ref{thm:growth} gives
the local gradient bound on a smaller ball about each stagnation
point. On that working ball,
Lemma~\ref{lem:variational-measure-perimeter} then gives a
nonnegative Radon measure and the weighted perimeter
estimate near stagnation points. None of these conclusions
requires an additional assumption on the free boundary.

For part~\textup{(i)}, these estimates allow us to apply the
monotonicity, compactness, and density results in
\cite[Theorem~3.5(ii), Lemmas~4.1--4.2 and Proposition~4.7]{VW11}.
The remaining issue is to exclude zero density.
Proposition~\ref{prop:no-cusps-variational} gives the two-dimensional
irrotational proof of \cite[Theorem~I]{DPY26}; it applies on the working ball with
the general gradient constant supplied on that ball by
Theorem~\ref{thm:growth}. It therefore replaces the cusp
argument in \cite[Lemma~4.4]{VW11}, which uses the strong
Bernstein estimate \eqref{eq.strong-bernstein}. This leaves precisely the Stokes and full-density
alternatives stated in part~\textup{(i)}.

For part~\textup{(ii)}, suppose that a full-density point exists.
The frequency result in \cite[Theorem~10.1]{VW11} gives the
normalized profile \eqref{eq:degenerate-profile}, which has an
interior nodal ray. The finite component assumption places nearby
zero points in the same connected zero component as the
stagnation point. Connectedness and the phase convergence in
part~\textup{(i)} force the interface to meet every circle in a
fixed range of radii above the line \(\{x_2=0\}\). Radial
projection then gives a positive lower bound for its
\(\mathcal H^1\) measure on a fixed disk. The singular set is
of zero \(\mathcal H^1\) measure, so the regular curve patches and their
one sided Bernoulli trace give a positive lower bound for the
mass of \(\Delta u\) on that disk. This contradicts \(u_r\to0\),
as proved in Proposition~\ref{prop:no-full-density}.
Once full density has been excluded, Weiss monotonicity and
convergence to the Stokes profile give local finiteness of the
stagnation set.

Part~\textup{(iii)} assumes the curve and trace conditions at
every free boundary point above the line \(\{x_2=0\}\).
These conditions are used in three places. They make the
extension by zero of \((|\nabla u_r|^2-x_2)^+\) continuous in
Lemma~\ref{lem:pressure-excess}; they allow continuation of
the interface arcs in Lemma~\ref{lem:arc-crossing}; and they
imply the viscosity boundary condition in
Lemma~\ref{lem:local-water-wave-flatness}. The resulting
estimates locate the interface near the Stokes rays and control
its normals uniformly on rescaled annuli. The final argument
then gives the two \(C^1\) branches.

\enlargethispage{3\baselineskip}

\begin{thebibliography}{99}

\bibitem{AC81}
H.~W.~Alt and L.~A.~Caffarelli,
\emph{Existence and regularity for a minimum problem with free boundary},
J. Reine Angew. Math. \textbf{325} (1981), 105--144.
\href{https://doi.org/10.1515/crll.1981.325.105}
{doi:10.1515/crll.1981.325.105}.

\bibitem{AFT82}
C.~J.~Amick, L.~E.~Fraenkel, and J.~F.~Toland,
\emph{On the Stokes conjecture for the wave of extreme form},
Acta Math. \textbf{148} (1982), 193--214.
\href{https://doi.org/10.1007/BF02392728}{doi:10.1007/BF02392728}.

\bibitem{CR13}
X.~Cabr\'e and X.~Ros-Oton,
\emph{Sobolev and isoperimetric inequalities with monomial weights},
J. Differential Equations \textbf{255} (2013), no.~11, 4312--4336.
\href{https://doi.org/10.1016/j.jde.2013.08.010}{doi:10.1016/j.jde.2013.08.010}.

\bibitem{CRS16}
X.~Cabr\'e, X.~Ros-Oton, and J.~Serra,
\emph{Sharp isoperimetric inequalities via the ABP method},
J. Eur. Math. Soc. \textbf{18} (2016), no.~12, 2971--2998.
\href{https://doi.org/10.4171/JEMS/659}{doi:10.4171/JEMS/659}.

\bibitem{DHN26}
G.~De Philippis, J.~Hirsch, and M.~Nahon,
\emph{Optimal regularity for vectorial, higher order and
non-minimizing Bernoulli problems},
preprint, arXiv:2609.12180v2, version~2 (29 September 2026).
\url{https://arxiv.org/abs/2609.12180v2}.

\bibitem{DS11}
D.~De Silva,
\emph{Free boundary regularity for a problem with right hand side},
Interfaces Free Bound. \textbf{13} (2011), no.~2, 223--238.
\href{https://doi.org/10.4171/IFB/255}{doi:10.4171/IFB/255}.

\bibitem{DPY26}
L.~Du, Y.~Pu, and C.~Yang,
\emph{Nonexistence of cusp singularities for gravity
waves of extreme form},
preprint, arXiv:2611.11776v1 (9 November 2026).
\url{https://arxiv.org/abs/2611.11776v1}.

\bibitem{EG15}
L.~C.~Evans and R.~F.~Gariepy,
\emph{Measure Theory and Fine Properties of Functions},
revised ed., Textbooks in Mathematics, CRC Press, Boca Raton, 2015.

\bibitem{Jost13}
J.~Jost,
\emph{Partial Differential Equations},
3rd ed., Graduate Texts in Mathematics 214,
Springer, New York, 2013.
\href{https://doi.org/10.1007/978-1-4614-4809-9}{doi:10.1007/978-1-4614-4809-9}.

\bibitem{KW24}
D.~Kriventsov and G.~S.~Weiss,
\emph{Rectifiability, finite Hausdorff measure, and compactness for
non-minimizing Bernoulli free boundaries},
Comm. Pure Appl. Math. \textbf{78} (2025), no.~3, 545--591.
First published online 13 October 2024.
\href{https://doi.org/10.1002/cpa.22226}{doi:10.1002/cpa.22226}.

\bibitem{McC25}
S.~McCurdy,
\emph{Non-existence of cusps for a free-boundary problem
for water waves},
Interfaces Free Bound. \textbf{27} (2025), no.~1, 1--11.
\href{https://doi.org/10.4171/IFB/527}{doi:10.4171/IFB/527}.

\bibitem{Nek21}
A.~I.~Nekrasov,
\emph{On steady waves},
Izv. Ivanovo-Voznesensk. Politekhn. Inst.
\textbf{3} (1921), 52--65 (in Russian).

\bibitem{P82}
P.~I.~Plotnikov,
\emph{Justification of the Stokes conjecture in the theory of
surface waves}, Dinamika Sploshn. Sredy \textbf{57}
(1982), 41--76 (Russian).
English translation: \emph{Proof of the Stokes conjecture in
the theory of surface waves}, Stud. Appl. Math.
\textbf{108} (2002), no.~2, 217--244.
\href{https://doi.org/10.1111/1467-9590.01408}
{doi:10.1111/1467-9590.01408}.

\bibitem{ST08}
E.~Shargorodsky and J.~F.~Toland,
\emph{Bernoulli free-boundary problems.}
Mem. Amer. Math. Soc. \textbf{196} (2008), no.~914.
\href{https://doi.org/10.1090/memo/0914}{doi:10.1090/memo/0914}.

\bibitem{St80}
G.~G.~Stokes,
\emph{Considerations relative to the greatest height of oscillatory
irrotational waves which can be propagated without change of form},
in \emph{Mathematical and Physical Papers}, Vol.~I,
Cambridge University Press, Cambridge, 1880, pp.~225--228.

\bibitem{VW11}
E.~V\u{a}rv\u{a}ruc\u{a} and G.~S.~Weiss,
\emph{A geometric approach to generalized Stokes conjectures},
Acta Math. \textbf{206} (2011), no.~2, 363--403.
\href{https://doi.org/10.1007/s11511-011-0066-y}
{doi:10.1007/s11511-011-0066-y}.

\end{thebibliography}
\end{document}